 \documentclass[reqno,11pt,a4paper]{amsart}

\usepackage{pgfplots}
\usepackage{caption}
\usepackage{subcaption}
\pgfplotsset{compat=1.12} 

\usepackage{amsfonts,amsmath,amssymb}
\usepackage[latin1]{inputenc}
\usepackage{dsfont}
\usepackage{enumerate}
\usepackage{xcolor}
\usepackage[all]{xy}

\def\notshow#1\notshowend{} %
\usepackage{hyphenat}

\usepackage{enumitem}
\usepackage{mathtools}

\usepackage{graphicx}

\usepackage{tikz}
\usepackage{tikz-3dplot}
\usepackage{pgfplots}

\usepackage{multirow}
\usepackage{array}

\newcommand{\C}{\mathcal{C}}

\usepackage{amssymb}
\usepackage{amsfonts}
\usepackage{mathrsfs}
\usepackage{hyperref}

\usepackage[normalem]{ulem} 

\def\bb#1\eb{\textcolor{blue}{#1}} 
\def\br#1\er{\textcolor{red}{#1}} %
\def\bm#1\em{\textcolor{purple}{#1}} %

\newcommand{\R}{\mathds R}

\newcommand{\gn}{|g|}

\newtheorem{thm}{Theorem}[section]

\theoremstyle{definition}

\newtheorem{rem}[thm]{Remark}

\newcommand{\ben}{\begin{enumerate}}
\newcommand{\een}{\end{enumerate}}
\newcommand{\bit}{\begin{itemize}}
\newcommand{\eit}{\end{itemize}}
\newcommand{\edoc}{\end{document}}

\title[Gielis superformula and wildfire models]{Gielis superformula and wildfire models 
}

\author[M. A. Javaloyes]{Miguel \'Angel Javaloyes}\address{Departamento de Matem\'aticas, \hfill\break\indent Universidad de Murcia, \hfill\break\indent Campus de Espinardo,\hfill\break\indent 30100 Espinardo, Murcia, Spain} \email{majava@um.es}

\author[E. Pend\'as-Recondo]{Enrique Pend\'as-Recondo}\address{Departamento de Matem\'aticas, \hfill\break\indent Universidad de Murcia, \hfill\break\indent Campus de Espinardo,\hfill\break\indent 30100 Espinardo, Murcia, Spain}\email{e.pendasrecondo@um.es}

\author[M. S\'anchez]{Miguel S\'anchez}\address{Departamento de Geometr\'{\i}a y Topolog\'{\i}a, Facultad de Ciencias, \hfill\break\indent Universidad de Granada,\hfill\break\indent Campus Fuentenueva s/n, \hfill\break\indent 18071 Granada, Spain}\email{sanchezm@ugr.es}

\begin{document}
\begin{abstract}
 Experimental data on the propagation of wildfires show that its short-time spread has a double semi-elliptical shape.
	 Our main goal is to show that this shape can be accurately approximated in polar coordinates by choosing suitable parameters in Gielis superformula and, then, implemented 
	  in a geometrical model  for wildfire propagation. In this model, 
	 the firefront can be
determined by computing the lightlike geodesics of a specific Finsler spacetime, and
we derive a concise and efficient expression of the geodesic equations in these coordinates. 
\end{abstract}

\maketitle



\noindent {\em Keywords}: {\rm 
wildfire propagation, Finsler metrics, polar coordinates, Gielis superformula, Finsler spacetimes.}

\vspace{-4mm}
\section{Introduction}
Models for predicting the evolution of a wildfire can be very useful in many different situations. If the model can be obtained faster than the spreading of the fire, the firefighters can use it to organize themselves in a more efficient way and, for example,  to evacuate the areas in danger. It is also useful to design prevention strategies, to choose better where to put  firewalls or to remove the excess vegetation which can increase the speed of the fire. The key point in these models is to make an accurate approximation of the speed of the fire in every direction depending on the weather conditions (moisture, temperature), vegetation, slope and especially the wind (see \cite{Rot}). All these factors make this speed highly anisotropic and classically, the shape for the subset of velocities has been  approximated by  an ellipse  having the major axis  in the direction of the wind (see \cite{R,R2,finney1998,tymstra2010}). On the other hand, in 1983, based on experimental data, H. E. Anderson suggested that the best shape for the fire velocities is a double semi-ellipse (see  \cite{An}). So it makes sense to look for families of curves with a similar shape to the double semi-ellipse that can be described by an explicit formula. A very good tool to generate these curves is Gielis Superformula \cite{gielis03}. We will see that when we choose certain parameters, the generated shape is quite similar to the experimental results obtained by Anderson.

Once we have a candidate for the accurate approximation of the double semi-ellipse, we have to compute the estimated firefront.  In \cite{M16,M17}, S. Markvorsen proposed to use Finsler geometry to compute  it  (rather than the  Richards' equations, developed in \cite{R}, typically used in wildfire simulators).   The idea is to choose a Minkowski norm in the surface $S$ where the fire is spreading, $F_p:T_pS\rightarrow [0,+\infty)$  at each point $p\in S$, such that the velocities of the fire have $F$-norm equal to one, that is, they provide the indicatrix of a Finsler metric. Then,  if the velocities  do not depend on time, the $F$-length of a curve $\gamma:[a,b]\rightarrow S$ provides the amount of time that the fire takes to travel this path $\gamma$. It follows that  geodesics departing orthogonally from the initial firefront at the instant $t_0$  and with $F$-length $t_1-t_0$ give  the firefront at the instant $t_1$ when the velocities of the fire do not vary on time. This is because the orthogonal geodesics to the firefront are those that minimize the distance, or in other words, the fastest paths that the fire follows. 
In order to consider time-dependent velocities, a non-relativistic spacetime $\R\times S$  endowed with the  Lorentz-Finsler metric
\[L(\tau,v)=\tau^2-F_t(v)^2\]
is introduced in \cite{JPS3}.
Here, the natural projection $t:\R\times S\rightarrow \R$ represents the absolute time. 
Given any curve $\gamma:[a,b]\rightarrow S$, if we want to describe how the fire travels along $\gamma$ we need to choose its reparametrization so that the  spacetime curve $t\rightarrow (t,\gamma(t))$  is lightlike, namely, $0=L(1,\dot \gamma)=1-F_t(\dot\gamma)^2$.  By the Finslerian relativistic Fermat principle (see \cite{Per06}), the  first-arriving  paths followed by the fire are lightlike pregeodesics, in this case, parametrized by the absolute time $t$ (see \cite[Cor. 4.6]{JPS3}).

In order to construct these models, at each instant $t$, one must specify the velocity of propagation of the fire at each point and direction, which makes it natural to use polar coordinates at each tangent space. Indeed,
it is easy to obtain
 the value of the associated Minkowski norm in polar coordinates $(r,\theta)$, namely, if the curve is described by $(v(\theta),\theta)$, then 
\[F(r,\theta)=\frac{r}{v(\theta)}.\]
As a matter of fact, Gielis superformula, which is also expressed in these coordinates, becomes a  very suitable  choice  to model fire velocities. 
From a technical point of view, however, it is natural to use polar coordinates only for velocities but not for the space itself. Indeed, to obtain a full change of coordinates betwen Cartesian and polar coordinates, one has to use expressions such as
\begin{align*}
r&=\sqrt{x^2+y^2},\\
\theta&=\arctan \left(\frac{y}{x}\right),
\end{align*}
and distinguish the cases in that $x>0$, $x<0$ and $x=0$,   or to use a more complicated formula for $\theta$, namely, 
\[\theta= 2  \arctan\left(
\frac{y}{x +\sqrt{x^2 + y^2}}\right),\]
valid in  $\R^2\setminus \{(x, 0) : x \leq  0\}$.  Taking into account the expression of Gielis superformula, any of these changes of variables will get very complicated. Our proposal is to make computations of all the Finsler quantities without applying the change of variables, remaining always in polar coordinates in the tangent space to represent the velocity, while keeping the Cartesian coordinates in the manifold to represent the position.

The paper is organized as follows. In \S \ref{Sect:prep}, we introduce some basic definitions about Finsler metrics and Finsler spacetimes, concluding with the geodesic equations \eqref{prelightlikegeo} for a class of Finsler spacetimes that it will be used to model wildfires. In \S \ref{Sect:polar}, we compute the fundamental tensor of a Minkowski norm determined by an indicatrix which is given in polar coordinates by $(v(\theta),\theta)$, for a certain function $v:[0,2\pi)\rightarrow (0,+\infty)$. Finally, we obtain the condition on the function $v(\theta)$ to generate a strongly convex curve, and therefore a Minkowski norm (see \eqref{eqvtheta} and \eqref{eq:u}). For the reader's convenience, in \S \ref{Sect:geoMin},  we deduce the geodesic equations \eqref{Eqgeo} of a Finsler surface $(S,F)$ (with no dependence on time) using polar coordinates in every tangent space. Therefore, these geodesic equations depend only on $v_p(\theta)$, the curve that determines the indicatrix at the tangent space $T_pS$ and its partial derivatives (with respect to $\theta$ and the coordinates of $p$). Next, in \S \ref{Sect:lightgeo}, we consider the time-dependent case $v^t_p(\theta)$, and then obtain the lightlike pregeodesic  equations \eqref{geolighttemp} parametrized by the absolute time $t$ of a certain non-relativistic Finsler spacetime $(\R\times S,L)$. Lastly, in \S \ref{Sect:Gielis}, we find some parameters of Gielis superformula (see \eqref{eq:superformula} and \eqref{eq:base}) that provide a double semi-elliptical shape. Moreover, in \S \ref{Sect:strongconv}, we  prove analytically  that these parameters generate a periodic strongly convex curve and therefore a Minkowski norm. Finally, we perform some simulations of wildfires with fire velocities modeled by Gielis superformula, both in the time-independent case (see Figure \ref{fig:simulation1}), using classical Finsler metrics, and in the time-dependent one, making use of a certain Finsler spacetime (see Figure \ref{fig:simulation2}).

\section{Preliminaries}\label{Sect:prep}
\subsection{Minkowski norms and Finsler metrics}
A Minkowski norm in a vector space $V$  is defined as a function $F:V\rightarrow [0,+\infty)$ such that 
\begin{enumerate}
\item $F$ is smooth away from the zero vector,
\item $F$ is positive homogeneous of degree $1$, namely, $F(\lambda v)=\lambda F(v)$ for all $\lambda>0$ and $v\in V$,
\item for every $v\in V\setminus \bf 0$, the fundamental tensor $g_v$ defined as
\begin{equation*} \label{fundten}
g_v(u,w)\coloneqq\frac{1}{2}\frac{\partial^2}{ \partial \eta\partial \delta} F(v+\eta u+\delta w)^2|_{\eta=\delta=0}
\end{equation*}
 $\forall u,w \in V$,  is a positive definite scalar product.
\end{enumerate}
Observe that the Minkowski norm is determined by its indicatrix
\begin{equation*}\label{indicatrix}
\Sigma:=\{v\in V: F(v)=1\}.
\end{equation*}
Moreover,  $1$ is a regular value of $F:V\rightarrow [0,+\infty)$, since by homogeneity (Euler's theorem), for $v\in \Sigma$,  $dF_v(v)=F(v)=1\not=0$. Moreover,  as the fundamental tensor restricted to the indicatrix is the second fundamental form of $\Sigma$ (see for example \cite[Eq. (2.3)]{JavSan11}), it turns out that a hypersurface $\Sigma$ is the indicatrix of a Minkowski norm if and only if it is diffeomorphic to the sphere, it contains the zero vector and it has positive curvature with respect to any Euclidean metric on $V$.

A Finsler metric in a manifold $M$ is defined as a function $F:TM\rightarrow [0,+\infty)$, where $TM$ is the tangent bundle of $M$, such that 
\begin{enumerate}
\item $F$ is smooth on $TM\setminus \bf 0$,
\item $F_p:=F|_{T_pM}:T_pM\rightarrow[0,+\infty)$ is a Minkowski norm for every $p\in M$.
\end{enumerate}
One of the main features of Finsler metrics are  their  geodesics, which can be defined as the critical points of the length functional
\[\ell_F(\gamma ):=\int_a^b F(\dot\gamma(s)) ds,\]
where $\gamma$ is a piecewise smooth curve $\gamma:[a,b]\rightarrow M$. It can be proved that geodesics of a Finsler metric are locally minimizers of the (non-symmetric)-distance
\[d_F(p,q):=\inf_{\gamma\in C_p^q(M,[a,b])} \ell_F(\gamma),\]
where $C_p^q(M,[a,b])$ is the space of  piecewise smooth curves $\gamma:[a,b]\rightarrow M$ with $\gamma(a)=p$ and $\gamma(b)=q$.  This means that for a small enough piece of geodesic $\gamma:[a,a+\varepsilon]\rightarrow M$, the length of $\gamma$ coincides with the distance $d_F(\gamma(a),\gamma(a+\varepsilon))$. Geodesics can also be obtained as auto-parallel curves of some classical connections (Chern, Cartan, Berwald, Hashiguchi), but for our purposes, it will be enough to obtain the geodesics using the formal Christoffel symbols. Indeed, if $(\Omega, \varphi)$ is a chart of $M$, we define the formal Christoffel symbols as
\begin{equation*}\label{formalChris}
\gamma^k_{ij}(x,y):=\frac{1}{2} g^{km}(x,y)\left(\frac{\partial g_{im}}{\partial x^j}(x,y)+\frac{\partial g_{mj}}{\partial x^i}(x,y)-\frac{\partial g_{ij}}{\partial x^m}(x,y)\right)
\end{equation*}
where $(x,y)=(x^1,\ldots,x^n,y^1,\ldots,y^n)$ are the coordinates of a vector $v\in TM$ in the natural coordinates  $(\Omega\times \R^n,\hat\varphi)$  determined by $(\Omega,\varphi)$ and 
\[g_{ij}(x,y):=g_{\hat\varphi^{-1}(x,y)}(\partial_i,\partial_j)\]
 are the coefficients of the fundamental tensor, being $\{g^{ij}\}_{i,j=1,\ldots,n}$ the inverse matrix of $\{g_{ij}\}_{i,j=1,\ldots,n}$.  In the last equation we have used the Einstein convention for summation of indices and we will  use  it in the following without warning.  If a curve $\gamma$ is expressed in coordinates as $(\gamma^1,\ldots,\gamma^n)$, then it is a geodesic if and only if
\begin{equation}\label{geoeq}
\ddot\gamma^k+\gamma^k_{ij}(\gamma,\dot\gamma)\dot\gamma^i\dot\gamma^j=0
\end{equation}
(see \cite[Eq. (5.3.2)]{BCS00}).
\subsection{Lorentz-Finsler metrics and Finsler spacetimes}
Let $M$ be a manifold, $TM$ its tangent bundle and $\pi:TM\rightarrow M$ the natural projection. Let  $A\subset TM\setminus 0$ be an open subset  such that $A\cap T_pM$ is
non-empty,
 connected, conic and salient (the latter two meaning that if $A$ contains $v$ then it also contains $\lambda v$ for $\lambda >0$ but not for $\lambda <0$) 
 for every $p\in M$ 
and  $\bar A \setminus 0$ is a smooth manifold with boundary $\mathcal C$ embedded in $TM\setminus 0$. 
We will say that $L:A\rightarrow (0,+\infty)$ is a (properly) {\em Lorentz-Finsler metric  } if
\begin{enumerate}
\item $L$ is positive homogeneous of degree $2$, namely, $L(\lambda v)=\lambda^2 L(v)$,
\item for every $v\in A$, the fundamental tensor defined as
\begin{equation*} \label{fundten2}
g^L_v(u,w):=\frac{1}{2}\frac{\partial^2}{\partial \eta\partial \delta} L(v+\eta u+\delta w)|_{\eta=\delta=0}
\end{equation*}
for all $u,w\in T_{\pi(v)}M$, has index $n-1$,
\item $L$ is smoothly extendible to the boundary $\mathcal C$ as zero and the fundamental tensor is nondegenerate along the extension.
\end{enumerate}
 We will say that $(M,L)$  is a \emph{Finsler spacetime} and $A$ is its timelike cone, while the boundary $\C$ is its lightlike cone.  We will need a very particular type of Finsler  spacetime  constructed in a product manifold $\R\times M$ using a family of time-dependent Finsler metrics on $M$. Let 
 $F:\R\times TM\rightarrow [0,+\infty)$ be a function which is smooth on $\R\times TM\setminus\bf 0$ and such that $F^t_p:T_pM\rightarrow [0,\infty)$, defined as $F^t_p(v)=F(t,v)$ for all $v\in T_pM$, is a Minkowski norm for all $t\in \R$ and $p\in M$. Then,  $L:\R\times \R\times TM\rightarrow \R$, defined as
 \begin{equation}\label{lorentz-finsler}
 L(t,\tau,v):=\tau^2-F^t_p(v)^2
 \end{equation}
 with $p=\pi(v)$, is a Lorentz-Finsler metric in  $A=\{ (t,\tau,v): \tau> F^t_p(v)\}$. Observe that $L$ does not necessarily satisfy the smoothness condition on $L$, as it will be smooth on the vectors $(t,\tau,0)$ if and only if $F^t_p$ comes from a Riemannian metric (see \cite{Warner65}). In any case, $L$ is smooth on lightlike vectors and on lightlike geodesics and it can be smoothened without changing the metric on a neighborhood of  the lightcone $\mathcal C$ (see \cite[Section 5.2]{JS20}).  Anyway, it will not be necessary to smoothen $L$ because, as  commented in the introduction, to compute the firefront at a certain instant, we will need to obtain the lightlike pregeodesics of \eqref{lorentz-finsler} parametrized with the absolute time $t$. It turns out that we have to introduce some formal Christoffel symbols for $F^t_p$, namely, with the same conventions for notation as in \eqref{formalChris},  in a chart $(\Omega,\varphi)$ of $M$, 
 \begin{equation*}
\gamma^k_{ij}(t,x,y):=\frac{1}{2} g^{km}(t,x,y)\left(\frac{\partial g_{im}}{\partial x^j}(t,x,y)+\frac{\partial g_{mj}}{\partial x^i}(t,x,y)-\frac{\partial g_{ij}}{\partial x^m}(t,x,y)\right),
\end{equation*}
 where $g_{ij}(t,x,y)$ are the coordinates of the fundamental tensor of the Min\-kows\-ki norm $F^t_p$  and $(x,y)$ are the coordinates of the natural chart $(\Omega\times\R^n,\hat\varphi)$ of $TM$.  Finally, the equations for lightlike pregeodesics $\gamma:[a,b]\rightarrow \R\times M$, $t\mapsto (t, \gamma^1(t),\ldots,\gamma^n(t))$ parametrized by $t$ are  obtained modifying suitably \eqref{geoeq} (see \cite[Th. 4.11]{JPS3}), namely, 
 \begin{equation}\label{prelightlikegeo}
 \ddot \gamma^k
  +\dot \gamma^i g^{jk}(t,\gamma,\dot\gamma)\frac{\partial g_{ij}}{\partial t}(t,\gamma,\dot\gamma) 
+\dot\gamma^i \dot\gamma^j\left(\gamma^k_{ij}(t,\gamma,\dot\gamma)
 +\frac{1}{2} \dot\gamma^k\frac{\partial g_{ij}}{\partial t}(t,\gamma,\dot\gamma)\right)=0.
 \end{equation}
\section{Polar coordinates}\label{Sect:polar}
Assume that we have a Minkowski norm $F$ on $\R^2$ and, in polar coordinates, the indicatrix of $F$ is given by $r=v(\theta)>0$.  It follows that  
\[F(r,\theta)=\frac{r}{v(\theta)}.\]
 We can  obtain the fundamental tensor of $F$ computing the Hessian of  $\frac{1}{2} F^2$ in polar coordinates. This means that 
 \begin{equation}\label{fundtenhess}
 g_{w}(X,Y)=\frac{1}{2}(X(Y(F^2))-\nabla_XY (F^2)),
 \end{equation}
 for $w\in \R^2$.
 Let us observe that as $(x^1,x^2)=(r\cos(\theta),r\sin(\theta))$, then 
 \begin{align*}
 \partial_r&= \cos(\theta)\partial_1+\sin(\theta)\partial_2,\\
  \partial_\theta&=r(-\sin(\theta)\partial_1+\cos(\theta)\partial_2),
 \end{align*}
  where $\partial_1,\partial_2$ denote the partial vector fields of the Cartesian coordinates and $\partial_r$ and $\partial_\theta$ those of polar coordinates. 
Moreover,    using that $\partial_1$ and $\partial_2$ are parallel, it follows that 
 \begin{align}
 \nabla_{\partial_r}\partial_r &=0,\label{eq1}\\
  \nabla_{\partial_r}\partial_\theta&=-\sin(\theta)\partial_x+\cos(\theta)\partial_y=\frac{1}{r}\partial_\theta,\label{eq2}\\
  \nabla_{\partial_\theta}\partial_\theta&=-r\cos(\theta)\partial_x-r\sin(\theta)\partial_y=-r\partial_r.\label{eq3}
 \end{align}
 Then, if $w=r(\cos(\theta)\partial_1+\sin(\theta)\partial_2)$ and using \eqref{fundtenhess} and \eqref{eq1},
 \begin{equation}\label{grr}
  g_{rr}(r,\theta):=g_{w}(\partial_r,\partial_r)=\frac{1}{2}\frac{\partial^2}{\partial r^2}\left(\frac{r^2}{v(\theta)^2}\right)=\frac{1}{v(\theta)^2},
 \end{equation}
 using \eqref{fundtenhess} and \eqref{eq2},
 \begin{align}
 g_{r\theta}(r,\theta)&:=g_{w}(\partial_r,\partial_\theta)=\frac{1}{2}\frac{\partial^2}{\partial r\partial\theta}\left(\frac{r^2}{v(\theta)^2}\right)-\frac{1}{2}\nabla_{\partial_r}\partial_\theta \left(\frac{r^2}{v(\theta)^2}\right)\nonumber\\
 &= r\frac{\partial}{\partial\theta}\left(\frac{1}{v(\theta)^2}\right)-\frac{r}{2}\frac{\partial}{\partial\theta} \left(\frac{1}{v(\theta)^2}\right)
 =\frac{r}{2}\frac{\partial}{\partial\theta} \left(\frac{1}{v(\theta)^2}\right)=-\frac{r\dot v(\theta)}{v(\theta)^3},\label{grtheta}
 \end{align}
 and using \eqref{fundtenhess} and \eqref{eq3},
 \begin{align}
 g_{\theta\theta}(r,\theta)&:=g_{w}(\partial_\theta,\partial_\theta)=\frac{1}{2}\frac{\partial^2}{\partial\theta^2}\left(\frac{r^2}{v(\theta)^2}\right)-\frac{1}{2}\nabla_{\partial_\theta}\partial_\theta \left(\frac{r^2}{v(\theta)^2}\right)\nonumber\\
&= \frac{r^2}{2}\frac{\partial^2}{\partial\theta^2}\left(\frac{1}{v(\theta)^2}\right)+\frac{r}{2} \frac{\partial}{\partial r}\left(\frac{r^2}{v(\theta)^2}\right)
 \nonumber\\&= \frac{r^2}{2}\left( 6 \frac{\dot v(\theta)^2}{v(\theta)^4}-2\frac{\ddot v(\theta)}{v(\theta)^3}+\frac{2}{v(\theta)^2}
 	\right)\nonumber\\
 	&=\frac{r^2}{v(\theta)^2}\left( 3 \frac{\dot v(\theta)^2}{v(\theta)^2}-\frac{\ddot v(\theta)}{v(\theta)}+1\right).\label{gthetatheta}
 \end{align}
When computing the fundamental tensor in a point of the indicatrix $(v(\theta),\theta)$ (in polar coordinates), we obtain
\begin{align*}
g_{rr}(v(\theta),\theta)&= \frac{1}{v(\theta)^2},\\
g_{r\theta}(v(\theta),\theta)&= -\frac{\dot v(\theta)}{v(\theta)^2},\\
g_{\theta\theta}(v(\theta),\theta)&= 3 \frac{\dot v(\theta)^2}{v(\theta)^2}-\frac{\ddot v(\theta)}{v(\theta)}+1.
\end{align*}
 Using the above values of the fundamental tensor,   we can compute the fundamental tensor in the tangent vector to the indicatrix $\dot v(\theta)\partial_r+\partial_\theta$:
\begin{align*}
g_{w}(\dot v(\theta)\partial_r+\partial_\theta,\dot v(\theta)\partial_r+\partial_\theta)&= \dot v(\theta)^2 g_{rr}(v(\theta),\theta)\nonumber\\&\quad +2 \dot v(\theta) g_{r\theta}(v(\theta),\theta)+g_{\theta\theta}(v(\theta),\theta)\nonumber\\
&=\frac{\dot v(\theta)^2}{v(\theta)^2}-2 \frac{\dot v(\theta)^2}{v(\theta)^2}+3 \frac{\dot v(\theta)^2}{v(\theta)^2}-\frac{\ddot v(\theta)}{v(\theta)}+1\nonumber
\\&= 2 \frac{\dot v(\theta)^2}{v(\theta)^2}-\frac{\ddot v(\theta)}{v(\theta)}+1.
\end{align*}
This implies that the function $v(\theta)$ generates a Minkowski norm in polar coordinates if and only if
\begin{equation}\label{eqvtheta}
2\dot v(\theta)^2-\ddot v(\theta)v(\theta)+v(\theta)^2>0,
\end{equation}
for all $\theta \in [0,2\pi)$ and $v$ is  smoohtly periodic in $[0,2\pi]$.
\begin{rem}
\label{rem:strongly_convex}
 Notice that 
   \eqref{eqvtheta} is equivalent to saying that the indicatrix $r=v(\theta)$ is strongly convex (i.e., its curvature is strictly positive everywhere), and it can be rewritten as
   \begin{equation}
   \label{eq:u}
   \ddot u(\theta)+u(\theta)>0,
   \end{equation}
  with $u(\theta):=1/v(\theta)$.
 \end{rem}
\section{Computing the geodesics of a Finsler metric  considering polar coordinates in the tangent spaces   }\label{Sect:geoMin}
Let $(S,F)$ be a $2$-dimensional Finsler manifold and $(\Omega,\varphi)$ a chart on $S$ with $\varphi=(x^1,x^2)$. 
Now we will consider polar coordinates in every tangent space of $S$, namely, 
we will transform the Cartesian coordinates $(w^1,w^2)$ of $w\in T_pS$, $p\in S$,  meaning that 
\[w=w^1\partial_1+w^2\partial_2,\]
where $\partial_1$ and $\partial_2$ are the partial vector fields associated, respectively, with the coordinates $x^1$, $x^2$ of the chart $(\Omega,\varphi)$. Observe that this change will not affect the coordinates of a point in $S$. As a consequence, if $\gamma:[a,b]\rightarrow \Omega\subset S$ is a geodesic of $(S,F)$ with image in the chart $(\Omega,\varphi)$, the coordinates $(\gamma^1,\gamma^2)$ will remain the same, but we will
 consider the functions $r:[a,b]\rightarrow [0,+\infty) $ and $\theta:[a,b]\rightarrow [0,2\pi)$ such that
\begin{align}\label{gamma1}
\dot\gamma^1=r\cos(\theta),\quad
\dot\gamma^2=r \sin(\theta).
\end{align}
Then
\begin{align*}
\ddot\gamma^1&=\dot r\cos(\theta)-r\dot\theta \sin(\theta) ,\\
\ddot\gamma^2&=\dot r \sin(\theta)+r\dot\theta \cos(\theta),
\end{align*}
and it follows that
\begin{align}\label{dotr}
\dot r=&\cos(\theta)\ddot\gamma^1+\sin(\theta)\ddot\gamma^2,\\
\dot\theta=&\frac{1}{r}(-\sin(\theta)\ddot\gamma^1+\cos(\theta)\ddot\gamma^2).\nonumber
\end{align}
Assume also that for every $p\in S$, the indicatrix of $F_p$ is given by $v_p(\theta)$ as a function in polar coordinates, namely, 
\[\Sigma_p=\{(v_p(\theta)(\cos(\theta)\partial_1+\sin(\theta)\partial_2): \theta\in [0,2\pi)\}\]
  and
 let us compute now $g_{ij}(w)=g_w(\partial_i,\partial_j)$, $w\in TS$, $i,j=1,2$, in polar coordinates. We have that 
 \begin{align*}
 \partial_1=&\cos(\theta)\partial_r-\frac{1}{r}\sin(\theta)\partial_\theta,\\
 \partial_2=& \sin(\theta)\partial_r+\frac{1}{r}\cos(\theta)\partial_\theta.
 \end{align*}
Then if  $w=r(\cos(\theta)\partial_1+\sin(\theta)\partial_2)$, using \eqref{grr}, \eqref{grtheta} and \eqref{gthetatheta}, it follows that 
 \begin{align*}
 g_{12}(p,r,\theta)&=g_{w}(\partial_1,\partial_2)\\&=\sin(\theta)\cos(\theta)(g_{rr}-\frac{1}{r^2} g_{\theta\theta})+\frac{1}{r}(\cos^2(\theta)-\sin^2(\theta))g_{r\theta}
 \\&= \frac{1}{v_p(\theta)^2}\left(-\frac{1}{2}\sin(2\theta)\left( 3 \frac{\dot v_p(\theta)^2}{v_p(\theta)^2}-\frac{\ddot v_{p}(\theta)}{v_p(\theta)}\right)
 -\cos(2\theta)\frac{\dot v_p(\theta)}{v_p(\theta)}\right),
 \end{align*}
 \begin{align*}
 g_{11}(p,r,\theta)&=g_{w}(\partial_1,\partial_1)\\&=
 \cos^2(\theta)g_{rr}-\frac{2}{r}\cos(\theta)\sin(\theta)g_{r\theta}+\frac{1}{r^2}\sin^2(\theta)g_{\theta\theta}
 \\&=
 \frac{\cos^2(\theta)}{v_p(\theta)^2}+\sin(2\theta)\frac{\dot v_p(\theta)}{v_p(\theta)^3}+\frac{\sin^2(\theta)}{v_p(\theta)^2}\left( 3 \frac{\dot v_p(\theta)^2}{v_p(\theta)^2}-\frac{\ddot v_{p}(\theta)}{v_p(\theta)}+1\right)\\
 &=\frac{1}{v_p(\theta)^2}
 \left(1 +\sin(2\theta)\frac{\dot v_p(\theta)}{v_p(\theta)}+  \sin^2(\theta)\left(3\frac{\dot v_p(\theta)^2}{v_p(\theta)^2}-\frac{\ddot v_{p}(\theta)}{v_p(\theta)}\right)\right),
 \end{align*}
 \begin{align*}
 g_{22}(p,r,\theta)&=g_{w}(\partial_2,\partial_2)\\&=
 \sin^2(\theta)g_{rr}+\frac{2}{r}\cos(\theta)\sin(\theta)g_{r\theta}+\frac{1}{r^2}\cos^2(\theta)g_{\theta\theta}
 \\&=
  \frac{\sin^2(\theta)}{v_p(\theta)^2}-\sin(2\theta)\frac{\dot v_p(\theta)}{v_p(\theta)^3}+\frac{\cos^2(\theta)}{v_p(\theta)^2}\left( 3 \frac{\dot v_p(\theta)^2}{v_p(\theta)^2}-\frac{\ddot v_{p}(\theta)}{v_p(\theta)}+1\right)\\
 &=\frac{1}{v_p(\theta)^2}\left(1 -\sin(2\theta)\frac{\dot v_p(\theta)}{v_p(\theta)}+  \cos^2(\theta)\left(3\frac{\dot v_p(\theta)^2}{v_p(\theta)^2}-\frac{\ddot v_{p}(\theta)}{v_p(\theta)}\right)\right).
 \end{align*}
 \subsection{The final expression of the geodesic equation}
 
  Let us introduce the functions 
 \[\varphi(p,\theta):=\frac{\dot v_p(\theta)}{v_p(\theta)},\quad 
  \phi(p,\theta):= \frac{\ddot v_{p}(\theta)}{v_p(\theta)},\]
    denoting $\varphi_{\cdot j}$ and $\phi_{\cdot j}$, respectively,  the partial derivatives of $\varphi$ and $\phi$ with respect to $x^j$ and 
   \[v_{\cdot j}(p,\theta):=\partial_j(v_\theta)(p),\]
   where $v_\theta(p)=v_p(\theta)$, for $j=1,2$. Moreover, denote 
 \[ f_i(\theta):=
 \begin{cases}
 \sin^2(\theta) & \text{if $i=1$,}\\
 \cos^2(\theta) & \text{if $i=2$}
 \end{cases}
 \]
 then
 \begin{align*}
  g_{12}(p,\theta)=&\frac{1}{v_p(\theta)^2}\left(-\varphi(p,\theta)\cos(2\theta)-\frac{1}{2}\sin(2\theta)( 3\varphi^2(p,\theta)-\phi(p,\theta))\right),
  \end{align*}
  \begin{align*}
  g_{ii}(p,\theta)=& \frac{1}{v_p(\theta)^2}\left(1 -(-1)^i \sin(2\theta)\varphi(p,\theta)+ f_i(\theta) (3\varphi^2(p,\theta)-\phi(p,\theta))\right),
 \end{align*}
 \begin{align*}
\frac{ \partial g_{12}}{\partial x^j}(p,\theta)=&\frac{-2v_{\cdot j}(p,\theta)}{v_p(\theta)} g_{12}(p,\theta)+\frac{1}{v_p(\theta)^2}\left(-\varphi_{\cdot j}(p,\theta)\cos(2\theta)\right.\nonumber\\&\quad\quad\quad\quad\quad\quad\left.- \frac{1}{2}\sin(2\theta)(6\varphi(p,\theta)\varphi_{\cdot j}(p,\theta)-\phi_{\cdot j}(p,\theta))\right),
 \end{align*}
 \begin{align*}
 \frac{ \partial g_{ii}}{\partial x^j}(p,\theta)=&\frac{-2v_{\cdot j}(p,\theta)}{v_p(\theta)} g_{ii}(p,\theta)\nonumber\\
 &+\frac{1}{v_p(\theta)^2}\left(-(-1)^i  \sin(2\theta)\varphi_{\cdot j}(p,\theta)+f_i(\theta)(
 6\varphi(p,\theta)\varphi_{\cdot j}(p,\theta)-\phi_{\cdot j}(p,\theta)
 )\right).
 \end{align*}
 Putting $\gn =g_{11}g_{22}-g_{12}^2$, the  formal  Christoffel symbols can be computed as 
 \begin{align*}
 \gamma^1_{11}(p,\theta)&=\frac{\gn}{2}\left(g_{22}\frac{\partial g_{11}}{\partial x^1}-g_{12}\left(2\frac{\partial g_{12}}{\partial x^1}-\frac{\partial g_{11}}{\partial x^2}\right)\right),\nonumber\\
 \gamma^1_{12}(p,\theta)&=\frac{\gn}{2}\left(g_{22}\frac{\partial g_{11}}{\partial x^2}
 -g_{12}\frac{\partial g_{22}}{\partial x^1}\right),\nonumber\\
 \gamma^1_{22}(p,\theta)&=\frac{\gn}{2}\left(g_{22}\left(2\frac{\partial g_{12}}{\partial x^2}-\frac{\partial g_{22}}{\partial x^1}\right)-g_{12}\frac{\partial g_{22}}{\partial x^2}\right),\nonumber\\
 \gamma^2_{22}(p,\theta)&=\frac{\gn}{2}\left(-g_{12}\left(2\frac{\partial g_{12}}{\partial x^2}-\frac{\partial g_{22}}{\partial x^1}\right)+g_{11}\frac{\partial g_{22}}{\partial x^2}\right),\nonumber\\
 \gamma^2_{12}(p,\theta)&=\frac{\gn}{2}\left(-g_{12}\frac{\partial g_{11}}{\partial x^2}+g_{11}\frac{\partial g_{22}}{\partial x^1}
 \right),\nonumber\\
 \gamma^2_{11}(p,\theta)&=\frac{\gn}{2}\left(-g_{12}\frac{\partial g_{11}}{\partial x^1}+g_{11}\left(2\frac{\partial g_{12}}{\partial x^1}-\frac{\partial g_{11}}{\partial x^2}\right)\right). 
 \end{align*}
 Finally, from \eqref{geoeq}, \eqref{gamma1} and \eqref{dotr}, it follows that the geodesics are obtained solving the system
 \begin{align*}
 \dot\gamma^1&=r\cos(\theta),\\
\dot\gamma^2&=r \sin(\theta),\\
 \dot r&=-\cos(\theta) \Gamma^1(\gamma^1,\gamma^2,\theta)-\sin(\theta) \Gamma^2(\gamma^1,\gamma^2,\theta)\\
 \dot\theta&=\frac{1}{r}\left( \sin(\theta) \Gamma^1(\gamma^1,\gamma^2,\theta)- \cos(\theta) \Gamma^2(\gamma^1,\gamma^2,\theta)\right),
 \end{align*}
 where 
 \begin{align*}
 \Gamma^i(p,\theta)=r^2(\cos^2(\theta)\gamma^i_{11}(p,\theta)+2\cos(\theta)\sin(\theta)\gamma^i_{12}(p,\theta)
 +\sin^2(\theta)\gamma^i_{22}(p,\theta))
 \end{align*}
 for $i=1,2$. Observe that we can omit the equation for $\dot r$, as we can consider unit geodesics and in that case $r=v_\gamma(\theta)$:
 \begin{align}
 \dot\gamma^1&=v_\gamma(\theta)\cos(\theta),\nonumber\\
 \dot\gamma^2&=v_\gamma(\theta) \sin(\theta),\nonumber\\
 \dot\theta&=\frac{1}{v(\theta)}\left( \sin(\theta) \Gamma^1(\gamma^1,\gamma^2,\theta) - \cos(\theta) \Gamma^2(\gamma^1,\gamma^2,\theta)\right),\nonumber\\
 r&=v_\gamma(\theta).\label{Eqgeo}
 \end{align}
 \section{Computing the lightlike geodesics of our Finsler spacetime}\label{Sect:lightgeo}
 Assume now that we have a family of time-dependent Finsler metrics $F^t$ on a $2$-dimensional manifold $S$. Then we can construct the Lorentz-Finsler metric $(\R\times S,L)$ as defined in \eqref{lorentz-finsler}. For our fire model purposes, we will need to compute its lightlike pregeodesics parametrized with $t$, and by \eqref{prelightlikegeo}, we will need to compute the fundamental tensor of $F^t$ and its formal Christoffel symbols.
 Most of the formulas of the last section remain valid in this case. We only have to take into account that, as we 
 have a family of Finsler metrics depending on $t\in\R$, then the corresponding indicatrix in polar coordinates depends also on $t$ and it is a function $v^t_p(\theta)$. This implies that the expressions for $g_{ij}$ of last section basically hold in this case. Indeed, let us define \[\varphi(t,p,\theta):=\frac{\dot v_p^t(\theta)}{v^t_p(\theta)},\quad  
 \phi(t,p,\theta):= \frac{\ddot v_{p}^t(\theta)}{v^t_p(\theta)},\]
and denote $\phi_{\cdot j}$, $\varphi_{\cdot j}$ the partial derivatives with respect to $x^j$ and $\phi_{\cdot t}$ and $\varphi_{\cdot t}$ the partial derivatives with respect to $t$. Moreover, 
 \[ v_{\cdot j}(t,p,\theta):=\partial_j(v_\theta)(t,p), \quad  v_{\cdot t}(t,p,\theta):=\partial_t(v_\theta)(t,p),\]
 where $\partial_t$ is the partial vector of $\R\times S$ with respect the first coordinate and $v_\theta(t,p):=v^t_p(\theta)$.
  \begin{align*}
 g_{12}(t,p,\theta)=&\frac{1}{v^t_p(\theta)^2}\left(-\varphi(t,p,\theta)\cos(2\theta)- \frac{1}{2}\sin(2\theta)( 3\varphi^2(t,p,\theta)-\phi(t,p,\theta))\right),
 \end{align*}
 \begin{align*}
 g_{ii}(t,p,\theta)=& \frac{1}{v^t_p(\theta)^2}\left(1 -(-1)^i \sin(2\theta)\varphi(t,p,\theta)+ f^i(\theta) (3\varphi^2(t,p,\theta)-\phi(t,p,\theta))\right),
 \end{align*}
 \begin{align*}
 \frac{ \partial g_{12}}{\partial x^j}(p,\theta)&=\frac{-2v_{\cdot j}(t,p,\theta)}{v^t_p(\theta)} g_{12}(t,p,\theta)
 +\frac{1}{v^t_p(\theta)^2}\left(-\varphi_{\cdot j}(t,p,\theta)\cos(2\theta)\right.\nonumber\\&\quad\quad\quad\left.- \frac{1}{2}\sin(2\theta)(6\varphi(t,p,\theta)\varphi_{\cdot j}(t,p,\theta,p)-\phi_{\cdot j}(t,p,\theta))\right),
 \end{align*}
 \begin{align*}
 \frac{ \partial g_{ii}}{\partial x^j}(t,p,\theta)=&\frac{-2v_{\cdot j}(t,p,\theta)}{v^t_p(\theta)} g_{ii}(t,p,\theta)
 +\frac{1}{v^t_p(\theta)^2}\left(-(-1)^i\sin(2\theta)\varphi_{\cdot j}(t,p,\theta)\right.\nonumber\\&\quad\quad\quad\left.+f_i(\theta)(
 6\varphi(t,p,\theta)\varphi_{\cdot j}(t,p,\theta)-\phi_{\cdot j}(t,p,\theta)
 )\right),
 \end{align*}
 \begin{align*}
\frac{ \partial g_{12}}{\partial t}(t,p,\theta)=&\frac{-2v_{\cdot t}(t,p,\theta)}{v^t_p(\theta)} g_{12}(t,p,\theta)
+\frac{1}{v^t_p(\theta)^2}\left(-\varphi_{\cdot t}(t,p,\theta)\cos(2\theta)\right.\nonumber\\&\left.-\frac{1}{2}\sin(2\theta) (6\varphi(\theta,p,t)\varphi_{\cdot t}(t,p,\theta)-\phi_{\cdot t}(t,p,\theta))\right),
 \end{align*}
 \begin{align*}
 \frac{ \partial g_{ii}}{\partial t}(t,p,\theta)=&\frac{-2v_{\cdot t}(t,p,\theta)}{v^t_p(\theta)} g_{ii}(t,p,\theta)
 +\frac{1}{v^t_p(\theta)^2}\left(- (-1)^i \sin(2\theta)\varphi_{\cdot t}(t,p,\theta)\right.\nonumber\\&\left.+f_i(\theta)(
 6\varphi(t,p,\theta)\varphi_{\cdot t}(t,p,\theta)-\phi_{\cdot t}(t,p,\theta)
 )\right).
 \end{align*}
 In an analogous way to the geodesic equations for Finsler metrics and using \eqref{dotr} and \eqref{prelightlikegeo}, we deduce that the lightlike pregeodesics are the solutions of
 \begin{align}
 \dot\gamma^1&=v_\gamma^t(\theta)\cos(\theta),\nonumber\\
 \dot\gamma^2&=v^t_\gamma(\theta) \sin(\theta),\nonumber\\
 \dot\theta&=\frac{1}{v^t_\gamma(\theta)}\left( \sin(\theta) \Gamma^1(t,\gamma^1,\gamma^2,\theta) - \cos(\theta) \Gamma^2(t,\gamma^1,\gamma^2,\theta)\right),\nonumber\\
 r&=v^t_\gamma(\theta),\label{geolighttemp}
  \end{align}
 where 
 \begin{align*}
 \Gamma^i(t,p,\theta)&=r^2(\cos^2(\theta)(\gamma^i_{11}(t,p,\theta)+\frac{1}{2} \dot\gamma^i(t)\frac{\partial g_{11}}{\partial t}(t,p,\theta))\\&\quad+2\cos(\theta)\sin(\theta)(\gamma^i_{12}(t,p,\theta)
 +\frac{1}{2} \dot\gamma^i(t)\frac{\partial g_{12}}{\partial t}(t,p,\theta))
 \\&\quad+\sin^2(\theta)(\gamma^i_{22}(t,p,\theta))+\frac{1}{2} \dot\gamma^i(t)\frac{\partial g_{22}}{\partial t}(t,p,\theta)))\\
 &\quad  +rg^{ij}\left(\cos(\theta)\frac{\partial g_{1j}}{\partial t}+\sin(\theta)\frac{\partial g_{2j}}{\partial t}\right)
 \end{align*}
 for $i=1,2$.  Here we have used that lightlike pregeodesics $\gamma$ parametrized with the universal time $t$ satisfy that $F^t(\dot\gamma)=1$ and then, as in the case of classical Finsler metrics, the equation for $r$ is not necessary. 
 
\section{Computations using Gielis superformula}\label{Sect:Gielis}
In this last section, we will apply all of the above to the case when the function $ v(\theta) $ is given by Gielis superformula, i.e.,
\begin{equation}
\label{eq:superformula}
v(\theta) = \lambda\left(\left\lvert\frac{\cos\left(\frac{m(\theta-\varphi)}{4}\right)}{a}\right\rvert^{n_2}+\left\lvert\frac{\sin\left(\frac{m(\theta-\varphi)}{4}\right)}{b}\right\rvert^{n_3}\right)^{-\frac{1}{n_1}},
\end{equation}
where $ a, b, n_1 \in \mathds{R}\setminus \{0\} $ and $ m, n_2, n_3 \in \mathds{R} $ are the usual parameters of the superformula, and we have added $ \lambda > 0 $ to control the size of the figure via a homothety, and $ \varphi \in [0,2\pi) $ to control its orientation via a rotation.

\subsection{Double semi-ellipse}
\begin{figure}
	\centering
	\includegraphics[width=0.9\textwidth]{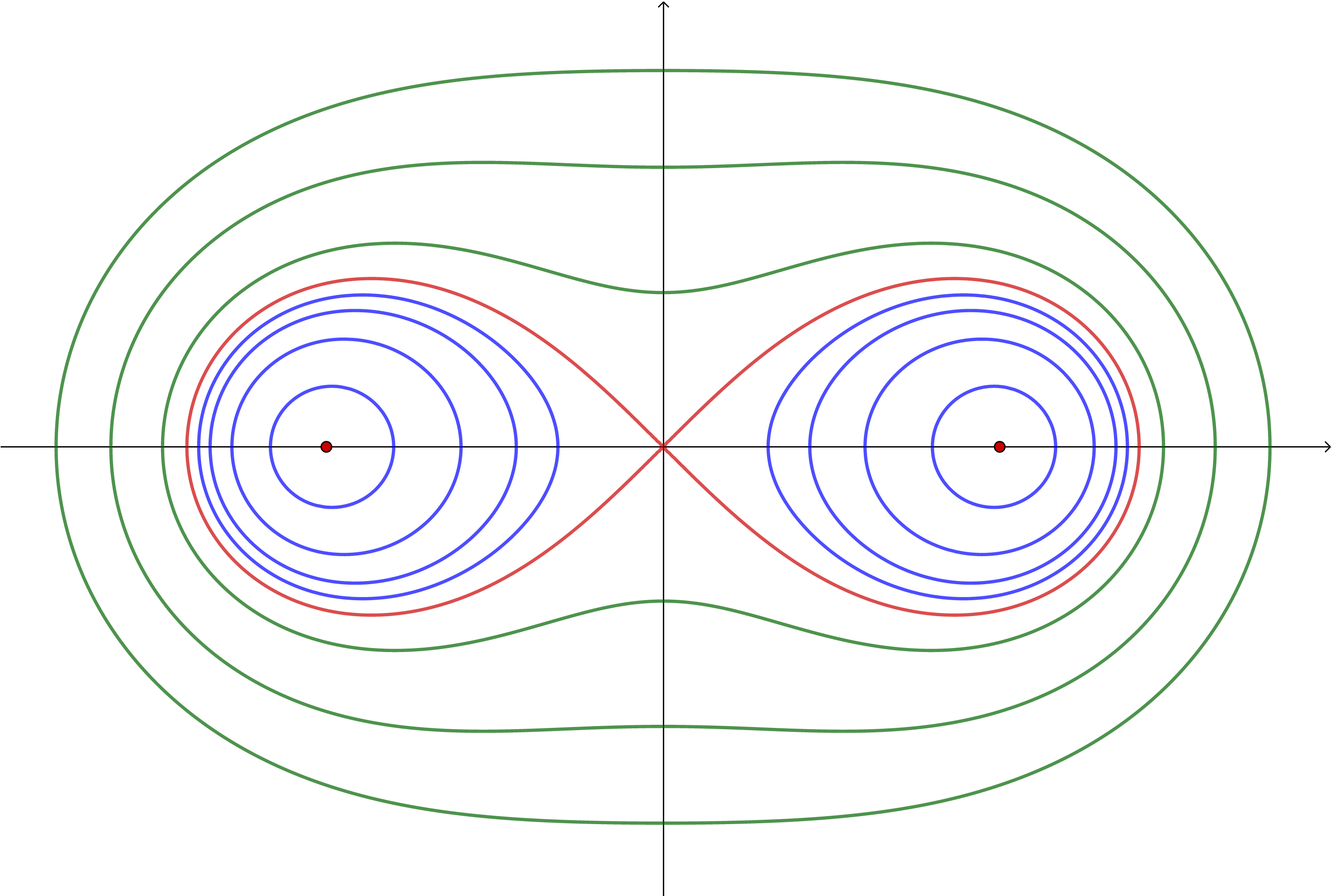}
	\caption{Cassini Ovals.}
	\label{fig:cassini}
\end{figure}

Gielis superformula features an incredibly wide range of complex shapes, depending on the choice of the parameters. In polar coordinates, these figures are given by the curve $ \theta \mapsto (v(\theta),\theta) $ and, among all the possibilities, we are especially interested in:
\begin{itemize}
\item Figures that are strongly convex, i.e., that satisfy \eqref{eqvtheta} or, equivalently, the condition \eqref{eq:u} in Remark~\ref{rem:strongly_convex}. This ensures that $ F(r,\theta) = r/v(\theta) $ is a proper Minkowski norm.
\item Figures resembling a double semi-ellipse,\footnote{Originally, the term \textit{double semi-ellipse} (or \textit{double ellipse}) refered specifically to the figure obtained by joining two semi-ellipses with a common minor axis. However, here we will be more lax in the use of the term (following \cite{JPS}), including similar shapes.} which is argued as the best experimental fitting for wildfire propagation in the presence of wind (see \cite{An}). Usually, the use of this shape complicates the computation of the firefront to the point that current fire growth simulators prefer to use a simple ellipse instead (see, e.g., \cite{finney1998,tymstra2010}). However, the innovative proposal, first introduced in \cite{M16}, of identifying the fire trajectories as Finsler geodesics---or, more generally, as lightlike pregeodesics of a Finsler spacetime, later proposed in \cite{JPS3}---, allows one to compute the wildfire propagation using formally the same ODE system---the geodesic equations---, regardless of the spread pattern (so long it is strongly convex). This obviously allows for more flexibility and accuracy, and new theoretical models incorporating the double semi-ellipse and other Finslerian patterns have been proposed recently (see, e.g., \cite{JPS}).
\end{itemize}
 In elementary differential geometry, 
 it is possible to find some classical curves with a certain resemblance to the double semi-ellipse  such as  
 some Cassini ovals. This is a family of quartics depending on two parameters, the semifocal distance $a$ and the constant product of the distances to the foci $b$, thus resembling the definition of the ellipse. The lemmniscata  (see the red curve in Fig. \ref{fig:cassini})  corresponds to the case of eccentricity $e (:=b/a)= 1$ and, for $e<1$, one has two Jordan curves. Each of them is strongly convex and (forgetting its partner curve) it exhibits a big asymmetry with respect to the two foci. However, the double ellipse effect is not too satisfactory (see  the blue curves in  Fig. \ref{fig:cassini})  and the possibilities to make the oval thiner to model a strong wind are more limited. This  suggests that we  look for  more suitable and adaptable curves using Gielis superformula.

\begin{figure}
	\centering
	\includegraphics[width=0.8\textwidth]{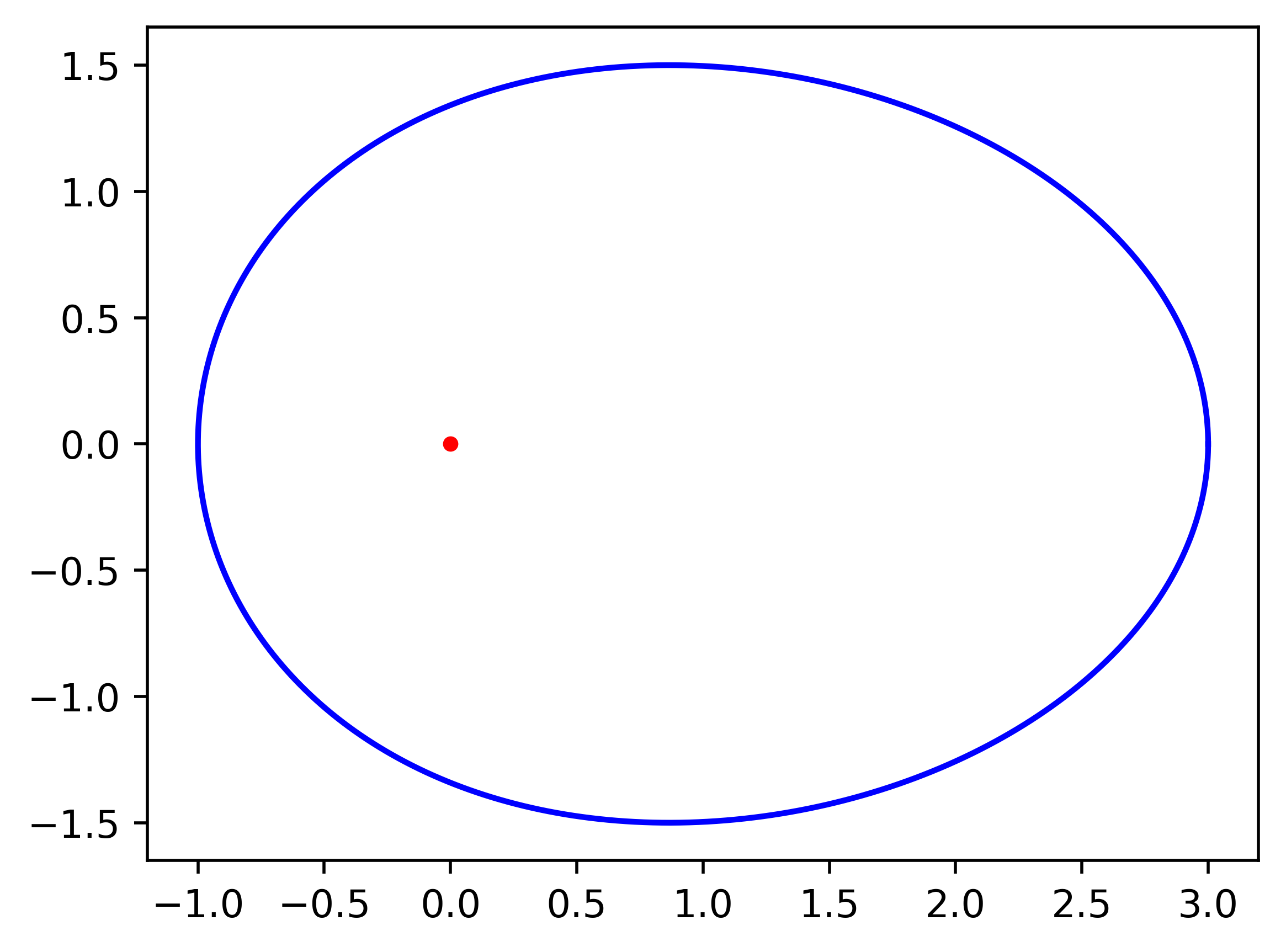}
	\caption{
Figure resembling a displaced ellipse given by Gielis superformula \eqref{eq:superformula} when $ a=3$, $ b=1$, $m=n_1=n_2=m_3=2$, $\lambda=1$ and $\varphi=0$.	}
	\label{fig:almost_semi-ellipse}
\end{figure}

 Choosing $m=4$ and $n_1=n_2=n_3=2$ in \eqref{eq:superformula}, one obtains an ellipse. However, the choice $m=2$ generates the asymmetry depicted in Fig. \ref{fig:almost_semi-ellipse} and, then, the choice $n_3=3$ yields the required double ellipse effect. 

\begin{figure}
\centering
\includegraphics[width=0.9\textwidth]{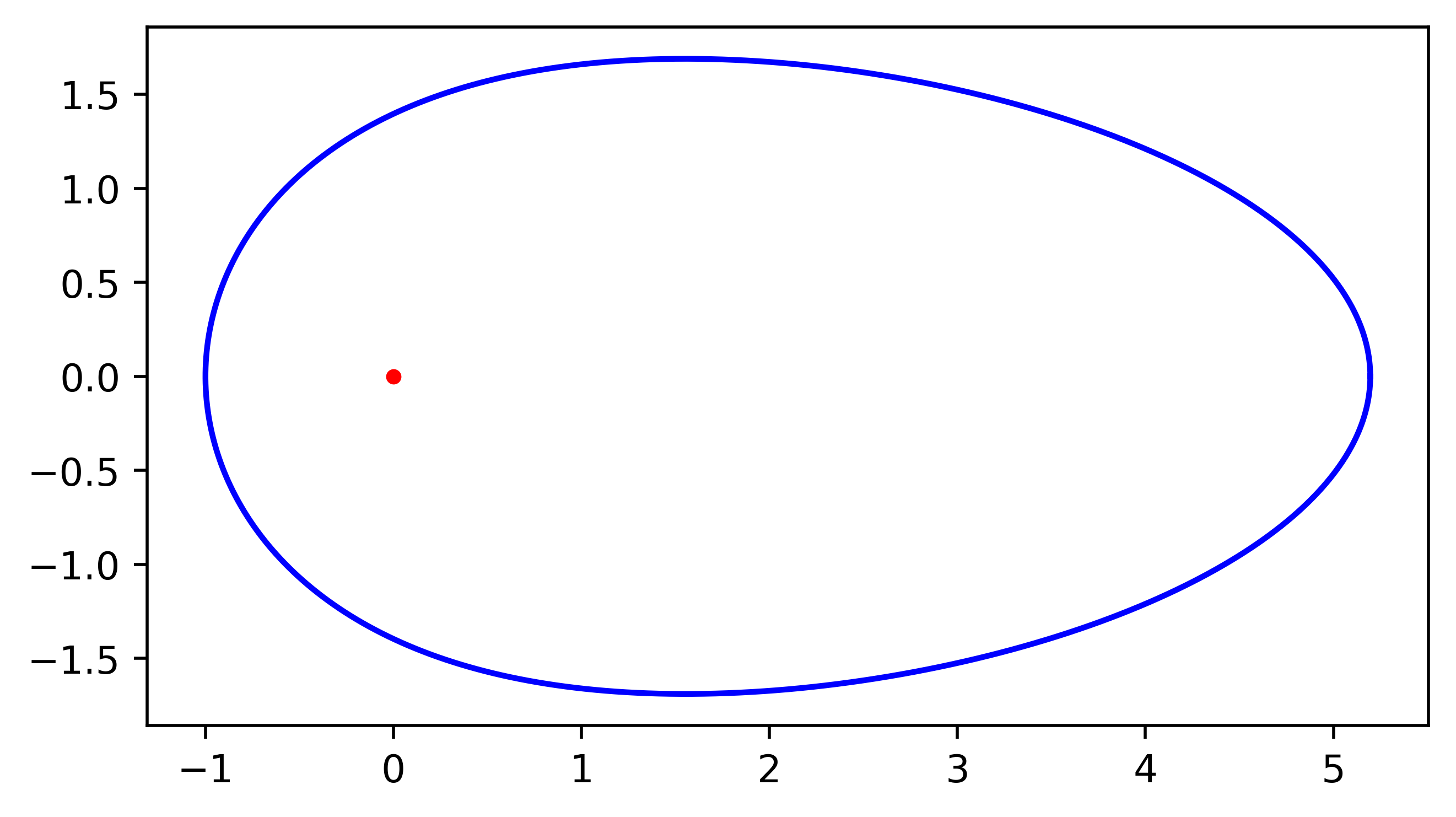}
\caption{Double semi-ellipse obtained by Gielis superformula \eqref{eq:superformula} when $ a = 3$, $b = 1$, $m = 2$, $n_1 = 2$, $n_2 = 3$, $n_3 = 2$, $\lambda = 1 $ and $ \varphi = 0 $.}
\label{fig:double_semi-ellipse}
\end{figure}

\begin{figure}
\centering
\begin{subfigure}{0.56\textwidth}
\includegraphics[width=\textwidth]{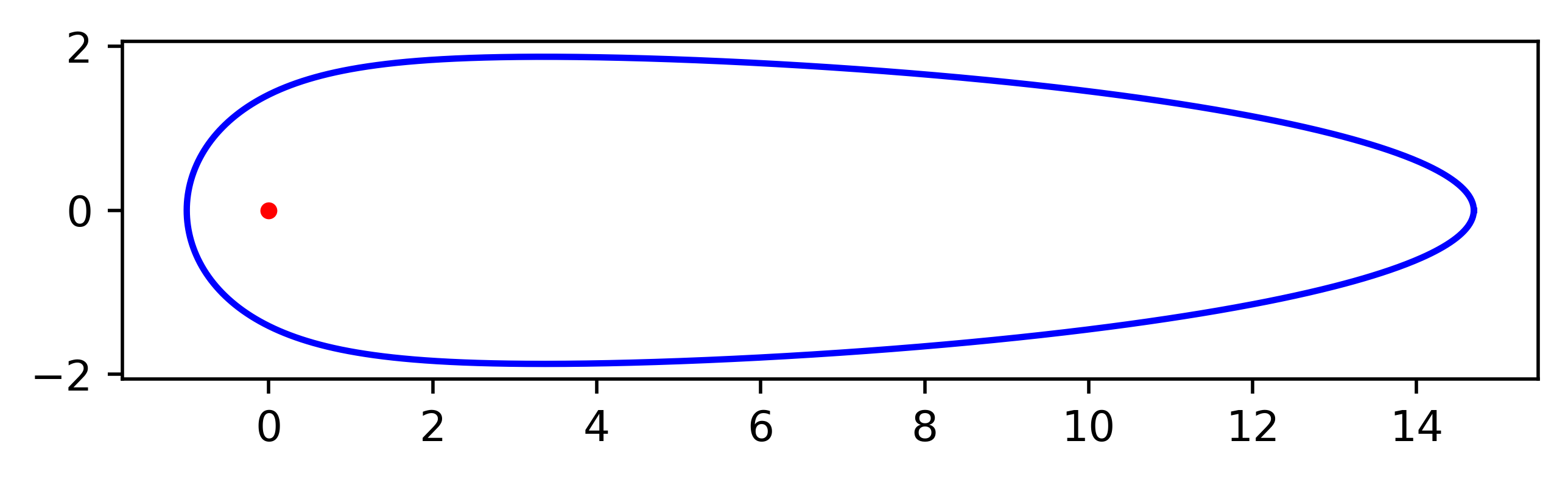}
\caption{$ a = 6 $.}
\label{fig:a6}
\end{subfigure}
\begin{subfigure}{0.42\textwidth}
\includegraphics[width=\textwidth]{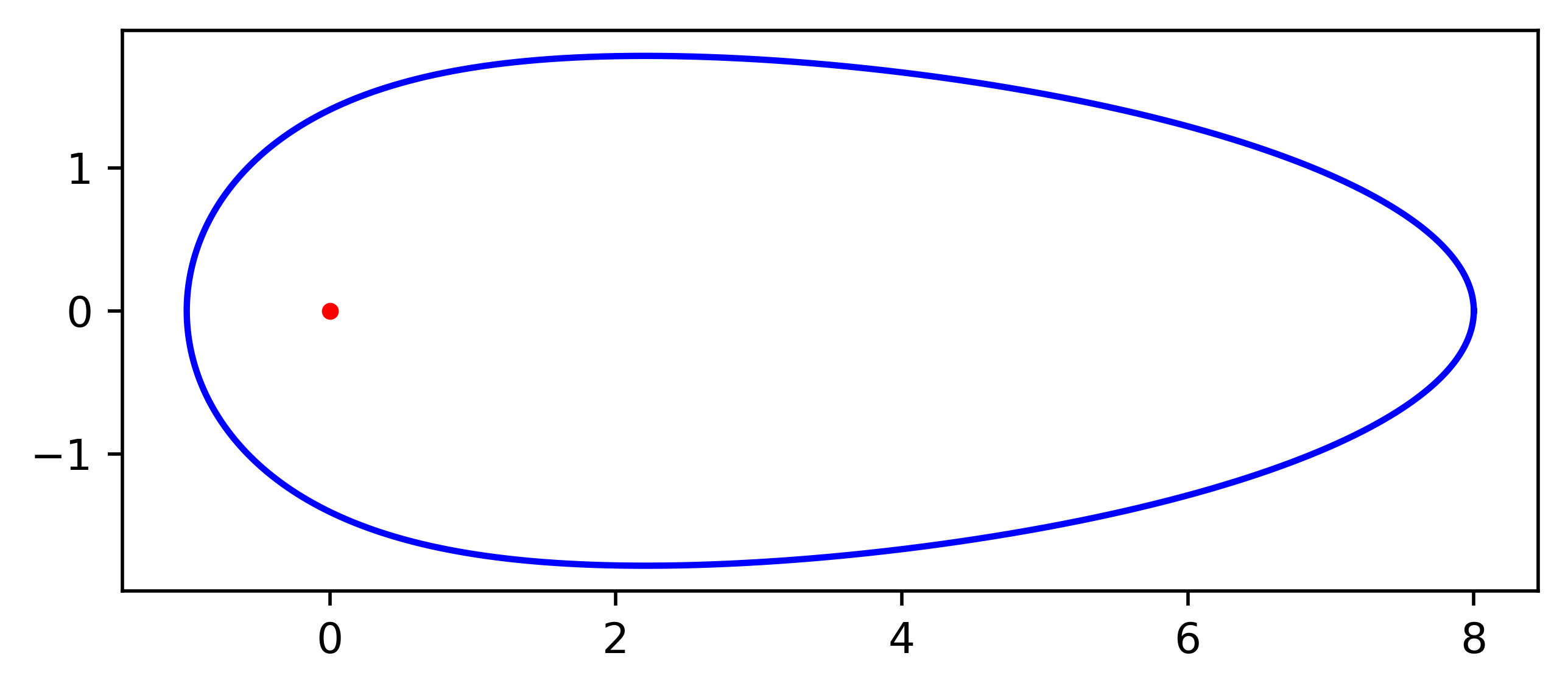}
\caption{$ a = 4 $.}
\label{fig:a4}
\end{subfigure}
\begin{subfigure}{0.56\textwidth}
\includegraphics[width=\textwidth]{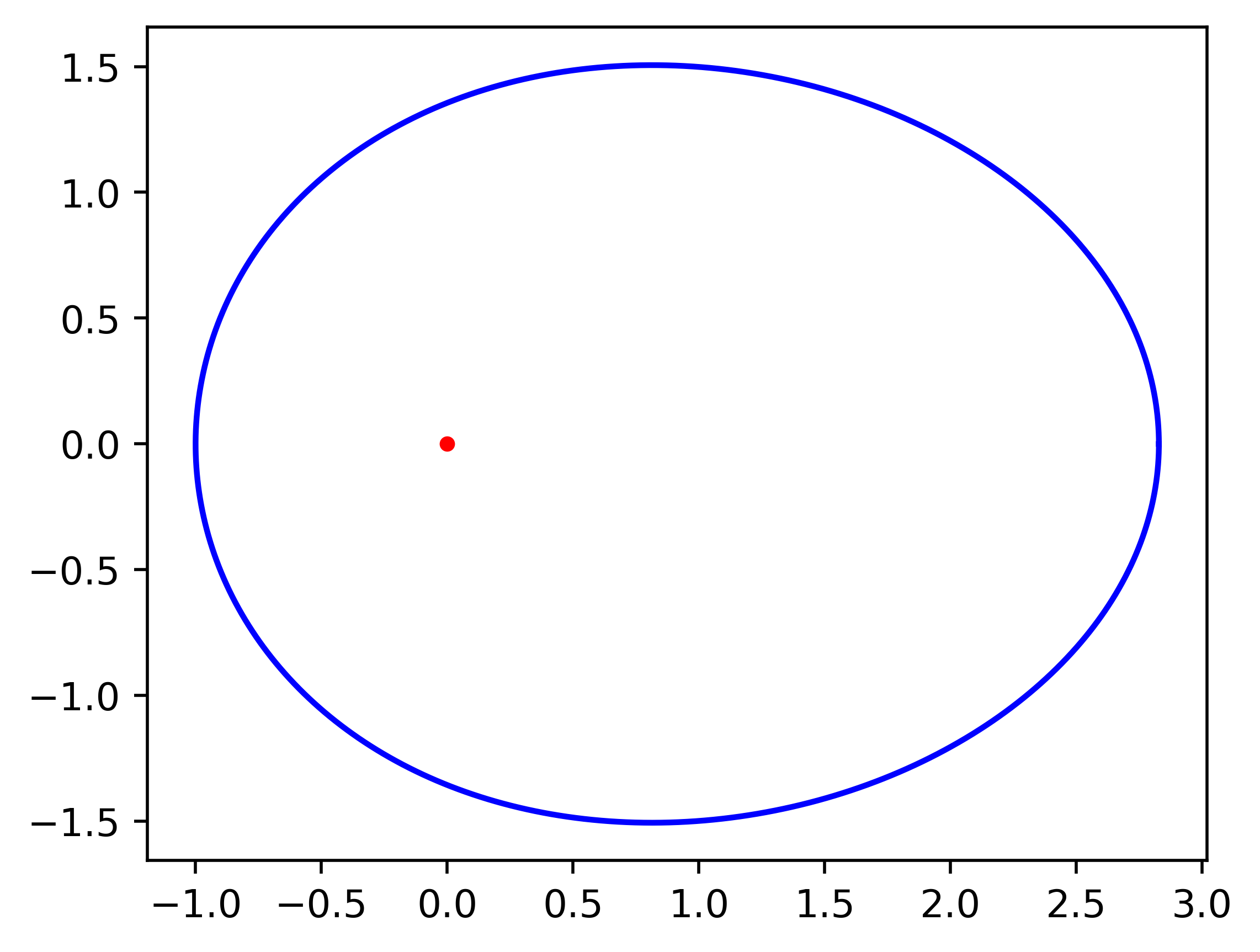}
\caption{$ a = 2 $.}
\label{fig:a2}
\end{subfigure}
\begin{subfigure}{0.42\textwidth}
\includegraphics[width=\textwidth]{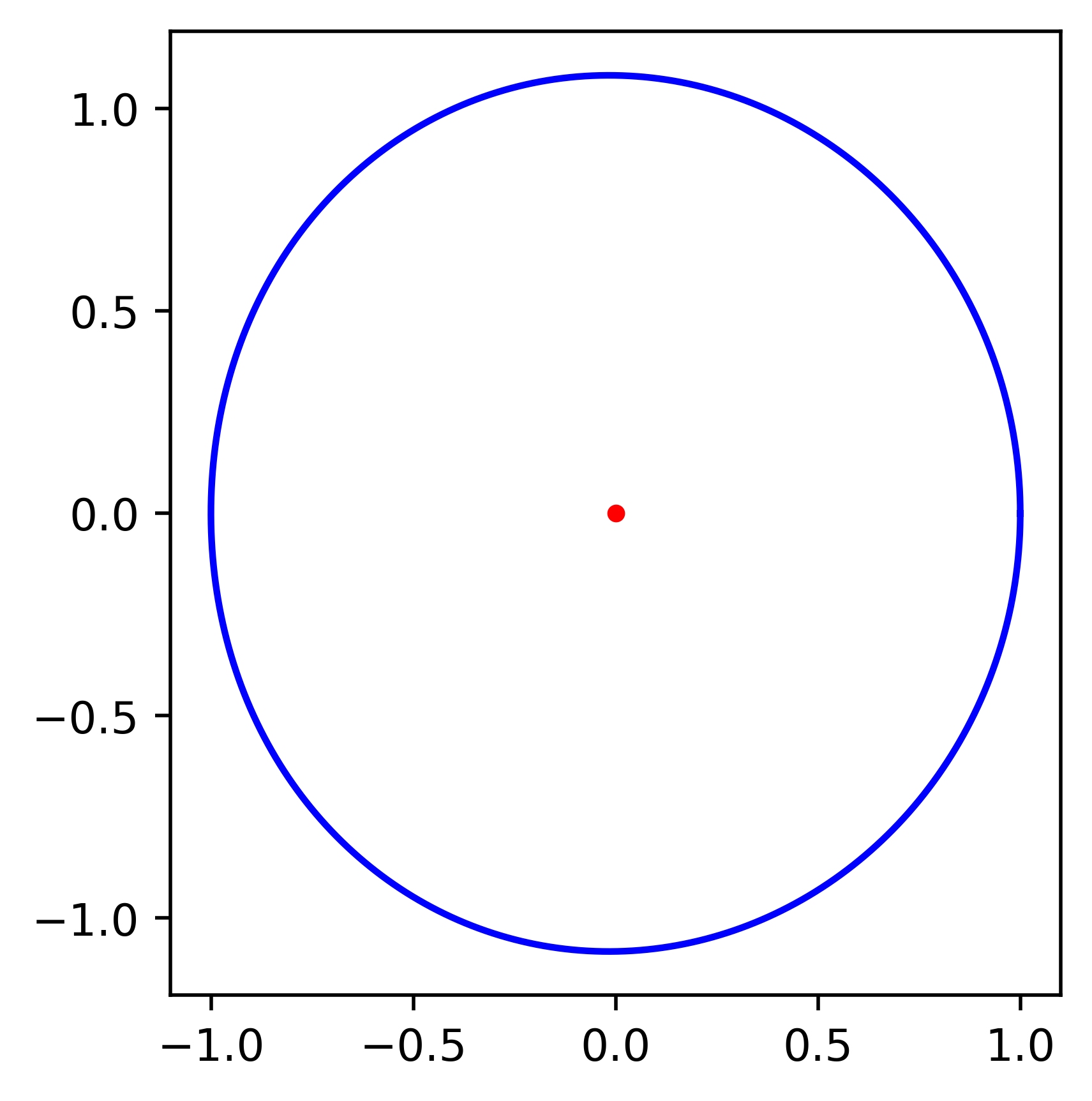}
\caption{$ a = 1 $.}
\label{fig:a1}
\end{subfigure}
\caption{Same choice of parameters as in Figure \ref{fig:double_semi-ellipse}, varying $ a $ from largest to lowest (with $ a \geq b $).}
\label{fig:a}
\end{figure}

\begin{figure}
\centering
\begin{subfigure}{0.44\textwidth}
\includegraphics[width=\textwidth]{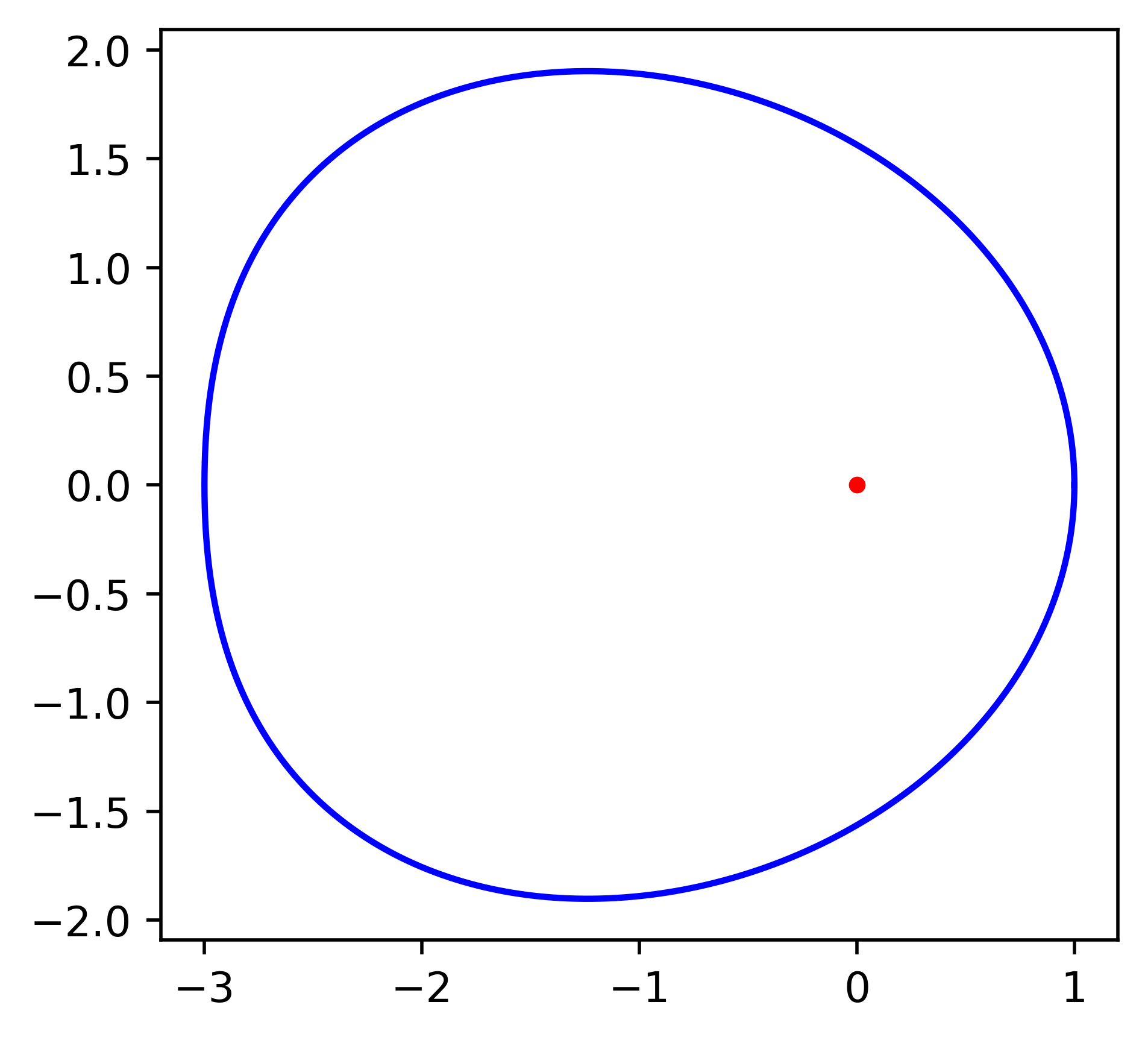}
\caption{$ b = 3 $.}
\label{fig:b3}
\end{subfigure}
\begin{subfigure}{0.54\textwidth}
\includegraphics[width=\textwidth]{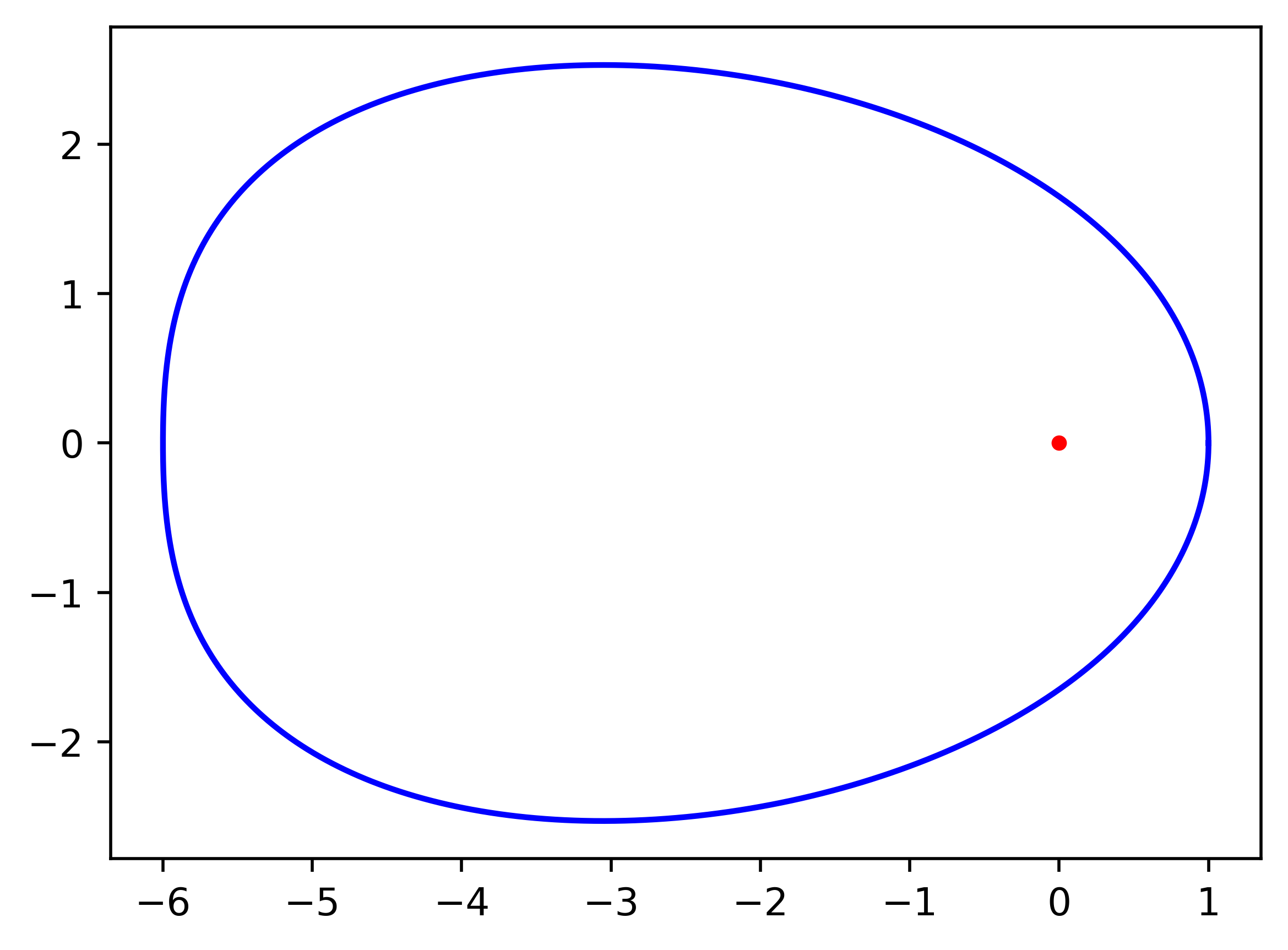}
\caption{$ b = 6 $.}
\label{fig:b6}
\end{subfigure}
\caption{Same choice of parameters as in Figure \ref{fig:double_semi-ellipse}, but now $ a = 1 $ and we vary the value of $ b $ (with $ b > a $).}
\label{fig:b}
\end{figure}

Keeping this in mind, it will be enough for our purposes here to fix the following values of parameters of $ v(\theta) $:
\begin{equation}
\label{eq:base}
m = 2, \quad n_1 = 2, \quad n_2 = 3, \quad n_3 = 2,
\end{equation}
and focus on the figures we obtain when varying $ a, b, \lambda $ and $ \varphi $. A good starting point is given by $ a = 3, b = 1, \lambda = 1 $ and $ \varphi = 0 $. The corresponding curve is shown in Figure~\ref{fig:double_semi-ellipse}, where one can observe its resemblance with a double semi-ellipse (compare with the double semi-elliptical fittings of real wildfires in \cite{An}, and with the theoretical double semi-ellipse constructed in \cite[Figure 4b]{JPS}). As mentioned before, the values of $ \lambda $ and $ \varphi $ simply change, respectively, the size of the figure and its orientation, keeping its shape. On the other hand, variations of $ a $ and $ b $ enable a wide range of different modifications to the shape of the figure. In principle $ a, b \in \mathds{R} \setminus \{0\} $, but notice that only their absolute values $ |a|, |b| $ play a role in \eqref{eq:superformula}, so we will assume $ a, b > 0 $, without loss of generality. So long $ a \geq b $, a general double semi-elliptical pattern is maintained: the larger the difference $ a-b $, the more elongated the double semi-ellipse (see Figures~\ref{fig:a6} and \ref{fig:a4}); the lower the difference, the more similar the figure becomes to a simple ellipse (see Figures~\ref{fig:a2} and \ref{fig:a1}). When $ b > a $, however, we abandon the double semi-ellipse to obtain a rather different shape (see Figure~\ref{fig:b}), which can nonetheless be also relevant for modeling other effects---apart from the wind---that take place in a wildfire, such as the slope (see the difference between wind-driven and slope-driven wildfires in \cite[Figure 5]{JPS}). Interestingly, notice that the direction of maximum speed can change from $ \theta = 0 $ when $ a > b $ (Figures~\ref{fig:a6}, \ref{fig:a4} and \ref{fig:a2}) to $ \theta = \pi $ when $ b > a $ (Figures~\ref{fig:b3} and \ref{fig:b6}; nevertheless, this does not always happen whenever $ b > a $), although this can be easily controlled by a suitable choice of $ \varphi $.

Obviously, the exact values of these parameters must be determined experimentally. In practice, the difference $ a-b $ should be related to the strength of the wind (possibly taking the slope into account, too), while the specific values of $ a, b $ and $ \lambda $ should depend on the fuel and meteorological conditions (composition of the soil, moisture, amount of oxygen in the air, etc.). Finally, $ \varphi $ should be given by a mixture between the wind direction
and the direction of maximum slope.


\subsection{Strong convexity}\label{Sect:strongconv}
Gielis superformula is known to provide a great variety of non-convex and even non-smooth figures, so it becomes crucial for our purposes here to ensure that we move within a range of parameters that always provides strongly convex shapes. The simplest way to check the strong convexity is to see whether condition \eqref{eq:u} holds. Although the complexity of $ v(\theta) $ in \eqref{eq:superformula} makes this condition very hard to analyze in general, it can be simplified by fixing some parameters of the superformula. Therefore, as stated above, we will restrict ourselves to the base values \eqref{eq:base}, varying only $ a, b \in \mathds{R} \setminus \{0\}, \lambda > 0 $ and $ \varphi \in [0,2\pi) $. At first glance, Figures~\ref{fig:double_semi-ellipse}, \ref{fig:a} and \ref{fig:b} seem to indicate that this choice of parameters might always yield a strongly convex shape, so let us verify this formally.

First, note that the values of $ \lambda $ and $ \varphi $ do not affect the strong convexity of the resulting figure, as this property is preserved under homothecies and rotations. Hence, for this analysis we can assume, without loss of generality, that $ \lambda = 1 $ and $ \varphi = 0 $. Putting $ u(\theta) \coloneqq 1/v(\theta) = \sqrt{h(\phi)} $, where $ v(\theta) $ is given by \eqref{eq:superformula} with fixed values \eqref{eq:base}, and
\begin{equation*}
h(\phi) := \frac{|\cos^3\phi|}{|a^3|} + \frac{\sin^2\phi}{b^2}, \quad \phi := \frac{\theta}{2},
\end{equation*}
then \eqref{eq:u} reduces to
\begin{equation*}
\begin{split}
\ddot{u}(\theta)+u(\theta) & = \frac{1}{8 \sqrt{h(\phi)}} \frac{d^2h(\phi)}{d\phi^2} - \frac{1}{16 h(\phi)^{3/2}}\left(\frac{dh(\phi)}{d\phi}\right)^2 + \sqrt{h(\phi)} \\
& = \frac{\frac{6|\cos\phi|\sin^2\phi}{|a^3|} - \frac{3|\cos^3\phi|}{|a^3|} - \frac{2\sin^2\phi}{b^2}  + \frac{2\cos^2\phi}{b^2}}{8\sqrt{\frac{|\cos^3\phi|}{|a^3|} + \frac{\sin^2\phi}{b^2}}} \\
& - \frac{\left(\frac{2\cos\phi \sin\phi}{b^2} - \frac{3|\cos\phi|\cos\phi\sin\phi}{|a^3|}\right)^2}{16\left(\frac{|\cos^3\phi|}{|a^3|} + \frac{\sin^2\phi}{b^2}\right)^{3/2}} + \sqrt{\frac{|\cos^3\phi|}{|a^3|} + \frac{\sin^2\phi}{b^2}} \\
& = \frac{\frac{3\sin^4\phi}{b^4} + \frac{1}{|a^3|b^2}\left(3|\cos\phi|\sin^4\phi + \frac{17}{2}|\cos^3\phi|\sin^2\phi + |\cos^5\phi|\right)}{4\left(\frac{|\cos^3\phi|}{|a^3|} + \frac{\sin^2\phi}{b^2}\right)^{3/2}} \\
& + \frac{\frac{1}{a^6}\left(\frac{3}{4}\cos^4\phi\sin^2\phi + \frac{5}{2}\cos^6\phi\right)}{4\left(\frac{|\cos^3\phi|}{|a^3|} + \frac{\sin^2\phi}{b^2}\right)^{3/2}}.
\end{split}
\end{equation*}
Observe that all the terms in the final expression are positive, so condition \eqref{eq:u} holds for all $ \theta \in [0,2\pi) $, i.e., the figure obtained by Gielis superformula, when restricted to the values \eqref{eq:base}, is strongly convex for all $ a, b \in \mathds{R} \setminus \{0\} $, $ \lambda > 0 $ and $ \varphi \in [0,2\pi) $.

\subsection{Simulations}
In the previous subsections we have seen that, taking Gielis superformula \eqref{eq:superformula} with fixed values \eqref{eq:base}, we can safely vary the paremeters $ a,b,\lambda $ and $ \varphi $ to obtain a wide variety of different double semi-ellipses (and even other Finslerian shapes), while keeping the strong convexity at all times. This ensures that
\begin{equation*}
F_p(r,\theta): = \frac{r}{v_p(\theta)}, \quad p \in S,
\end{equation*}
with $ v_p(\theta) = v(\theta,a(p),b(p),\lambda(p),\varphi(p)) $ given by \eqref{eq:superformula}, defines a proper Finsler metric on the 2-dimensional manifold $ M $, and if we include a time dependence, then
\begin{equation*}
L_{(t,p)}(\tau,r,\theta) := \tau^2 - F_p^t(r,\theta)^2, \quad (t,p) \in \mathds{R} \times S,
\end{equation*}
with $ F_p^t = r/v_p^t(\theta) $ and $ v_p^t(\theta) = v(\theta,a(t,p),b(t,p),\lambda(t,p),\varphi(t,p)) $ given by \eqref{eq:superformula}, defines a  Lorentz-Finsler metric on the 3-dimensional spacetime $ \mathds{R} \times S$ (recall the comments below \eqref{lorentz-finsler}).

When $ v_p^t(\theta) $ effectively represents the speed of the fire at every point $ p \in S $, time $ t \geq 0 $ and direction $ \theta \in [0,2\pi) $, then the curve $ \Sigma_p^t: \theta \mapsto (v_p^t(\theta),\theta) $, which is the indicatrix of $ F_p^t $, becomes the infinitesimal spread of the wildfire---i.e., $ \Sigma_p^t $ would be the firefront after one time unit if the parameters $ a(t,p),b(t,p),\lambda(t,p),\varphi(t,p) $ were independent of $ p $ and $ t $ (and in this case, the fire trajectories would be straight lines). When the fire speed depends only on the position, the fire trajectories can be identified as geodesics of $ F$ (see \cite{M16}), which can be obtained via the ODE system \eqref{Eqgeo}. When there is, in addition, a time-dependence, the fire trajectories become the projection on $ S $ of $ t $-parametrized lightlike pregeodesics of $ L $ (see \cite{JPS3,JPS}), which can be obtained via the ODE system \eqref{geolighttemp} (obviously, both ODE systems coincide in the time-independent case).\footnote{Here we are focusing specifically on wildfire propagation, but this can be applied in general to any wave.}

When performing wildfire simulations, one solves the ODE system for a finite amount of time $ t_0 > 0 $ and a finite number of fire trajectories, departing from a fixed ignition point $ p_0 \in S $ (or more generally, from a submanifold acting as an initial firefront, in which case the fire trajectories must depart $F$-orthogonally to it). If there are no intersections between the fire trajectories, then the last point of each trajectory becomes a point in the firefront at $ t_0 $, which can then be interpolated as a smooth curve---of course, the more trajectories we compute, the more precision we will have in the end to reconstruct the firefront. If there are intersections, however, it means that one or more trajectories have lost the property of being time-minimizing at some point (i.e., beyond this point, they only travel through already burned regions because other fire trajectories have arrived earlier) and therefore, they must not be taken into account to obtain the final firefront.\footnote{Technically, a fire trajectory remains in the firefront (i.e., remains time-minimizing) so long it does not reach its {\em cut point} (see \cite[\S 4]{JPS3}), although in practice one looks for intersections between fire trajectories (see \cite[\S 5.2]{JPS}).}

\begin{figure}
\centering
\includegraphics[width=0.9\textwidth]{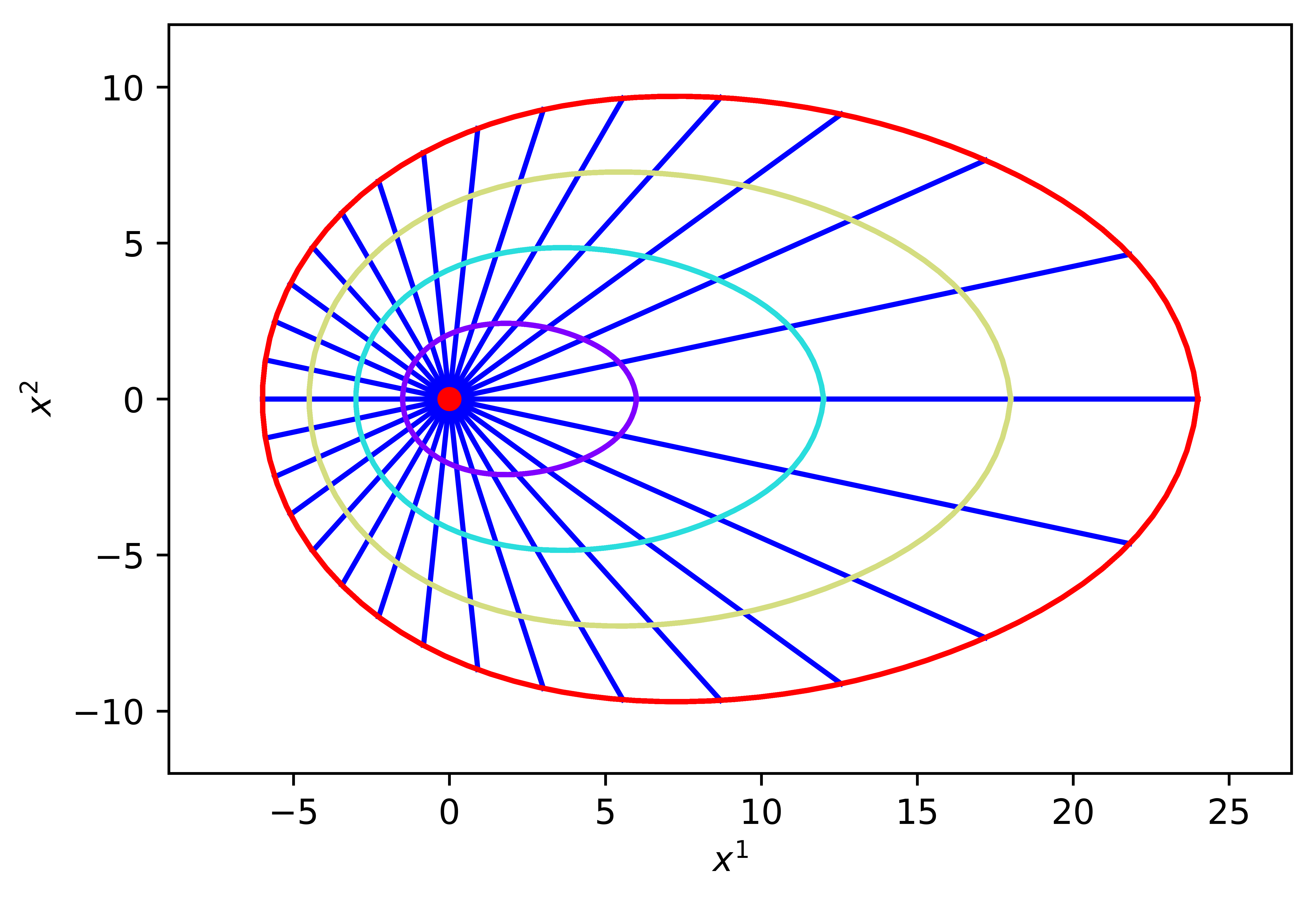}
\caption{Parameters: $ a = 4$, $b = 2$, $\lambda = 1$, $\varphi = 0$. Ignition point: $ p_0 = (0,0) $. Time: $ t_0 = 3 $.}
\label{fig:simulation0}
\end{figure}

\begin{figure}
\centering
\includegraphics[width=0.9\textwidth]{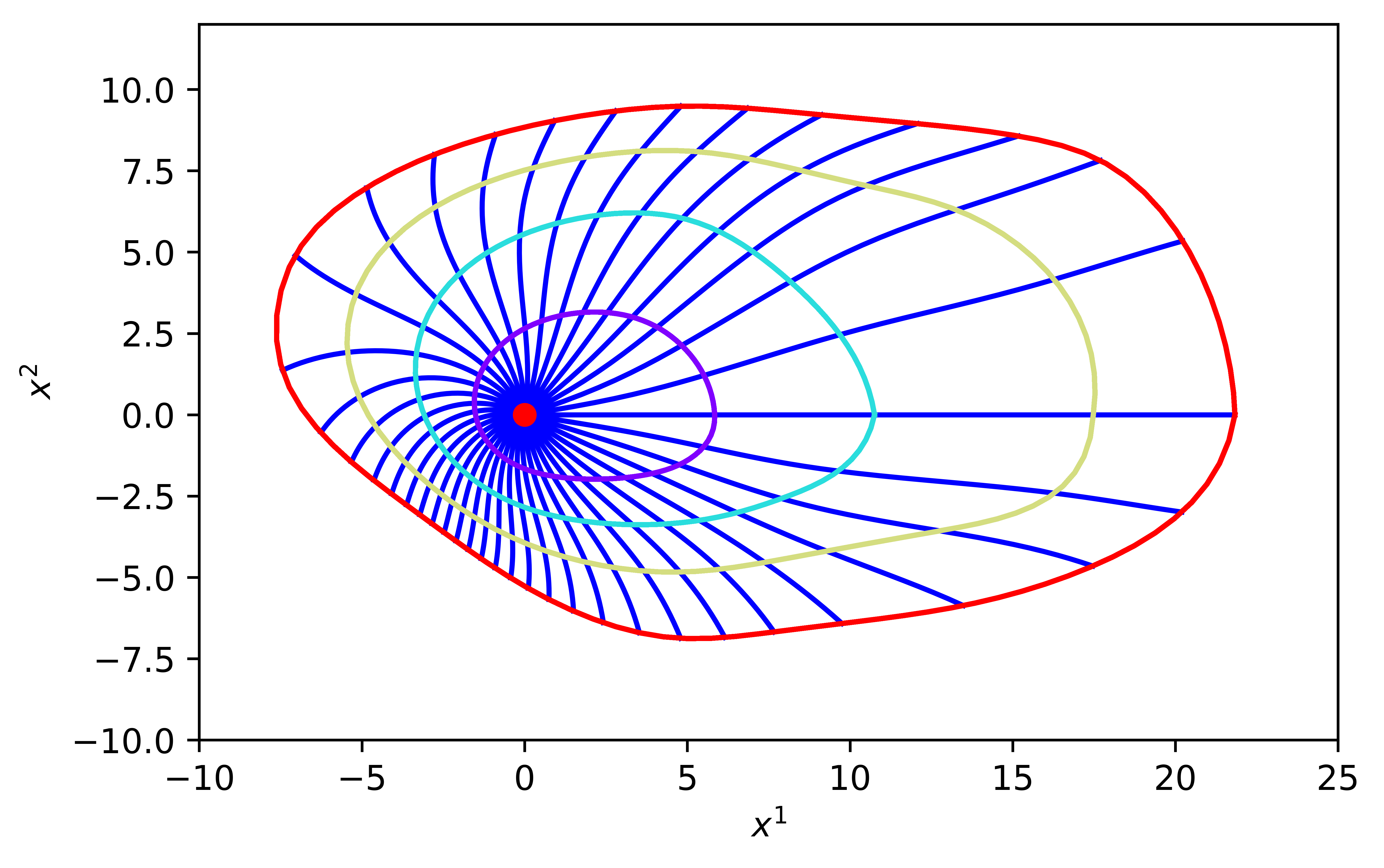}
\caption{Parameters: $ a = 4+\cos(x^1/2)$, $b = 2+\sin(x^2/2)$, $\lambda = 1$, $\varphi = 0$. Ignition point: $ p_0 = (0,0) $. Time: $ t_0 = 3 $.}
\label{fig:simulation1}
\end{figure}

\begin{figure}
\centering
\begin{subfigure}{0.49\textwidth}
\includegraphics[width=\textwidth]{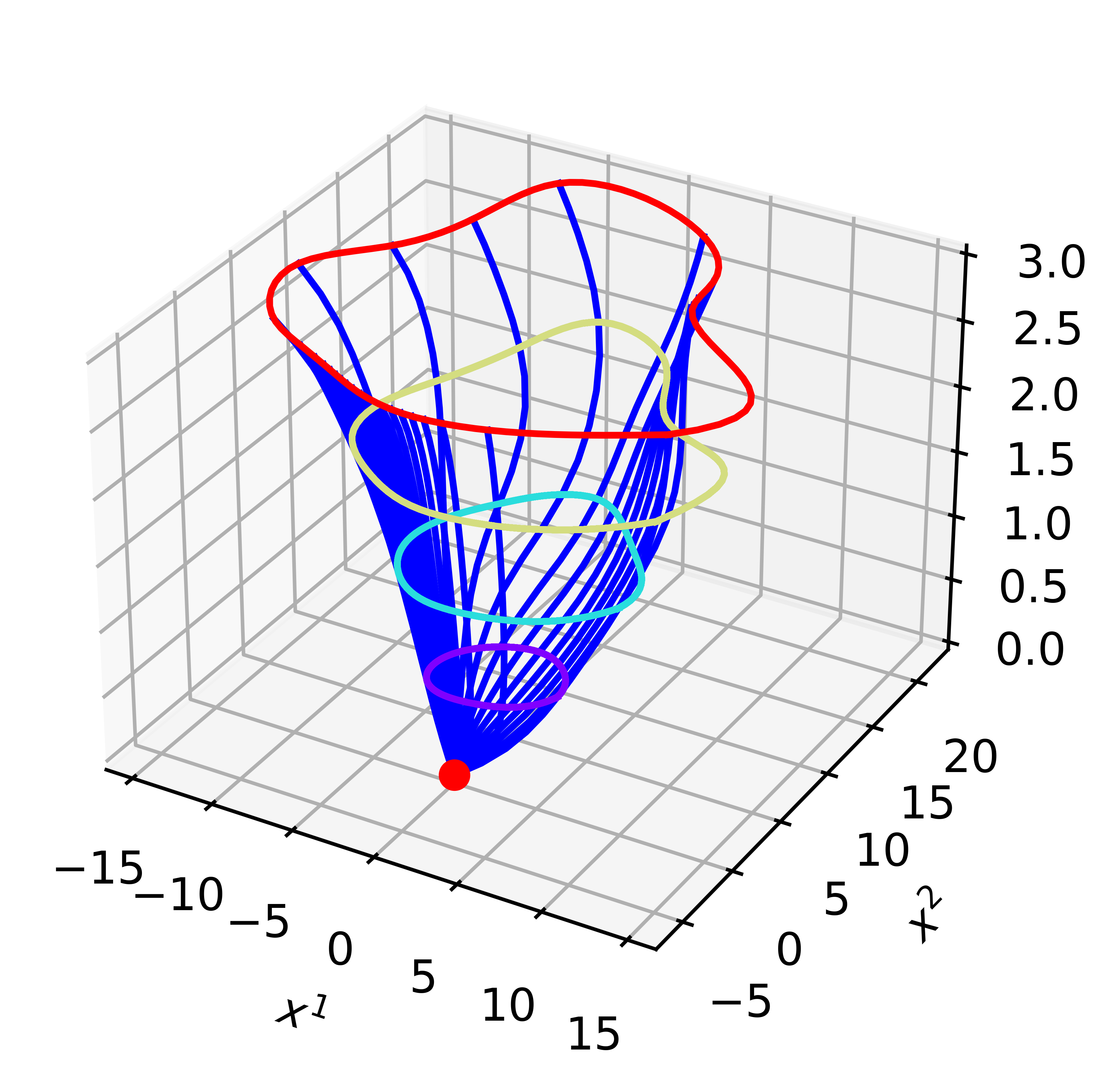}
\caption{Spacetime trajectories.}
\label{fig:sim_t2}
\end{subfigure}
\begin{subfigure}{0.5\textwidth}
\includegraphics[width=\textwidth]{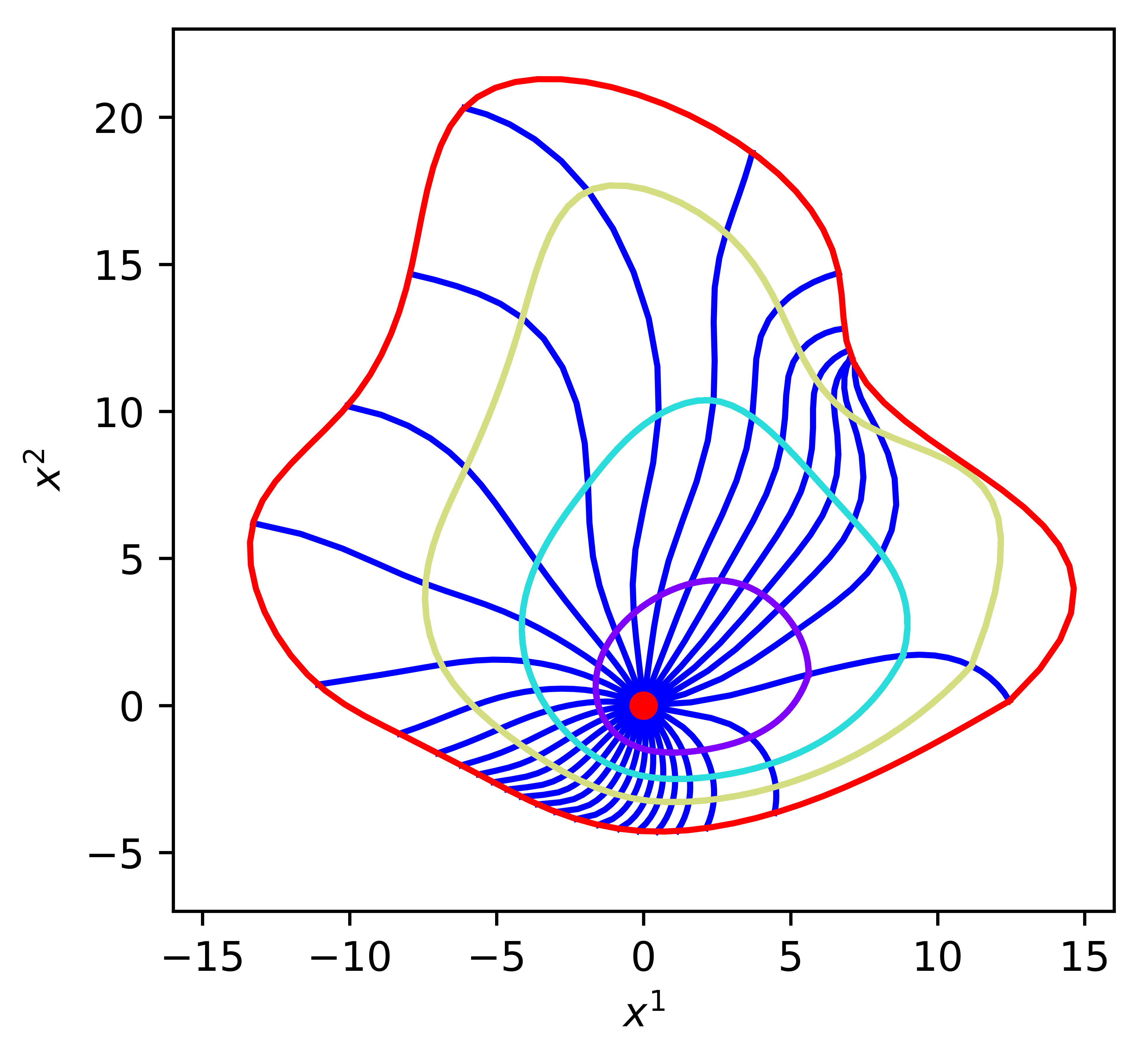}
\caption{Projection on the space.}
\label{fig:sim2}
\end{subfigure}
\caption{Parameters: $ a = 4+\cos(x^1/2)+t/2$, $b = 2+\sin(x^2/2)$, $\lambda = 1$, $\varphi = t$.  Ignition point: $ p_0 = (0,0) $. Time: $ t_0 = 3 $.}
\label{fig:simulation2}
\end{figure}

Figures~\ref{fig:simulation0}, \ref{fig:simulation1} and \ref{fig:simulation2} display some theoretical examples of these simulations, from least to greatest complexity, showing the firefront at regular intervals of time. In Figure~\ref{fig:simulation0} every parameter is constant, so the fire trajectories are straight lines and the firefronts coincide with the indicatrix scaled up. In Figure~\ref{fig:simulation1}, $ a $ and $ b $ depend on the position, and in Figure~\ref{fig:simulation2} we include an additional time-dependence in $ \varphi $. In this last situation, the actual solutions to the ODE system are lightlike pregeodesics in the spacetime (see Figure~\ref{fig:sim_t2}), which can then be projected on the space (see Figure \ref{fig:sim2}). Notice how, with each step, the fire trajectories become more complex and less similar to straight lines, so intersections are more likely to show up at later times (they cannot appear in the first instants of time by \cite[Theorem 4.8]{JPS3}).

\section*{Acknowledgments}
MAJ and EPR were partially supported by the project PID2021-124157NB-I00, funded by MCIN/AEI/10.13039/501100011033/``ERDF A way of making Europe", and also by Ayudas a proyectos para el desarrollo de investigaci\'{o}n cient\'{i}fica y t\'{e}cnica por grupos competitivos (Comunidad Aut\'{o}noma de la Regi\'{o}n de Murcia), included in the Programa Regional de Fomento de la Investigaci\'{o}n Cient\'{i}fica y T\'{e}cnica (Plan de Actuaci\'{o}n 2022) of the Fundaci\'{o}n S\'{e}neca-Agencia de Ciencia y Tecnolog\'{i}a de la Regi\'{o}n de Murcia, REF. 21899/PI/22. EPR and MS were partially supported by the project PID2020-116126GB-I00 funded by MCIN/AEI/10.13039/501100011033  as well as the framework IMAG-Mar\'{i}a de Maeztu grant CEX2020-001105-M/AEI/10.13039/501100011033. EPR was also supported by Ayudas para la Formaci\'{o}n de Profesorado Universitario (FPU) from the Spanish Government.

\end{document}